\documentclass[preprint,12pt]{elsarticle}

\usepackage{relsize}
\usepackage{hyperref}
\hypersetup{
    colorlinks=true,
    linkcolor=blue,
    filecolor=magenta,
    urlcolor=cyan,
}

\usepackage{amsmath}
\usepackage{amssymb}
\usepackage{setspace}
\usepackage[section]{placeins}
\usepackage{titlesec}
\usepackage{subcaption}
\titleformat{\paragraph}
{\normalfont\normalsize}{\theparagraph}{1em}{}
\titlespacing*{\paragraph}
{0pt}{3.25ex plus 1ex minus .2ex}{1.5ex plus .2ex}
\usepackage{amsthm}
\usepackage[left=3cm,right=2cm,top=2.5cm,bottom=2cm]{geometry}
\usepackage{epstopdf}
\journal{Numerical Methods for Partial Differential Equations}

\begin{document}

\begin{frontmatter}



\title{Adomian Decomposition Based Numerical Scheme for Flow simulations}

\author[bcam]{Imanol Garcia-Beristain~\corref{cor1}}
\ead{igarcia@bcamath.org}
\cortext[cor1]{Corresponding author: Tel.:+34 946 567 842 }

\author[bcam,alfaisal]{Lakhdar Remaki}

\address[bcam]{BCAM - Basque Center for Applied Mathematics, Alameda Mazarredo 14, 48009, Bilbao, Spain}

\address[alfaisal]{Department of Mathematics and Computer Science, Alfaisal University, KSA}

\begin{abstract}
This paper proposes a numerical method based on the Adomian decomposition
approach for the time discretization, applied to Euler equations. A recursive
property is demonstrated that allows to formulate the method in an appropriate
and efficient way. To obtain a fully numerical scheme, the space discretization
is achieved using the classical DG techniques. The efficiency of the obtained
numerical scheme is demonstrated through numerical tests by comparison to exact
solution and the popular Runge-Kutta DG method results.
\end{abstract}


\end{frontmatter}



\section{Introduction}

Linearized Euler equations (LEE) are extensively used in many problems
modeling and simulation. In particular simulating wave propagation in
aeroacoustic field. These equations offer, on one hand, the advantage to be
faster than solving the nonlinear Euler equations since large domains are
required for propagation. On the other hand, the information about the mean
flow is preserved comparing to a simple wave equation. This is due to the fact
that the linearization is achieved around a mean flow. Many aeroacoustic
applications are of big importance for industry and for human life quality
improvement. For instance, noise reduction in transportation. Especially, with
the sensitive population's mobility growth thanks to the development of fast
transportation facilities. In the US, the Joint Planning and Development Office
(JPDO) is planning a new NextGen system that increases the air traffic by a
factor of 3 towards 2025. Therefore reducing harmful sound effects becomes
critical in order to achieve this goal because the JPDO indicates that without
this substantial effort the number of people exposed to very high levels of
noise will increase substantially.

The LEE equations are numerically solved using different methods
including finite volume, finite elements, discontinuous Galerkin and spectral
methods. Research is still actively ongoing to design more effective numerical
schemes suitable for the large scales required by the practical problems.

In this paper we propose an accurate and cost-effective numerical scheme based
on the semi-analytical Adomian decomposition method, proposed by Adomian
\cite{adm1, adm2, adm3}, and used by many authors to solve a big range of
problems ranging from linear or nonlinear equations for deterministic or
stochastic PDEs.
It seems to be a promising trend in the field of PDEs solution
approximation. See for instance the important work of Wazwaz \cite{adm4, adm5,
admbook1, adm6, adm7}, and other relevant articles \cite{adm8, adm9, adm10}.
The reader is especially referred to the nice review paper on the topic
\cite{adm5}.

To derive the proposed approach, the semi-analytical Adomian decomposition
technique is applied to the Euler equations. Then, a recursive
property for the obtained time scheme is proved. This allows to formulate the
method in a simple and practical form, easy to implement and cost-effective. To
fully derive the numerical scheme (including space discretization) the classical
discontinuous Galerkin (DG) approach is used. We refer to the obtained time
scheme by ABS standing for Adomian Based Schemes and ABS-DG when the space
discretization is achieved using discontinuous Galerkin method. To demonstrate
the ABS-DG effectiveness, the (DG) method is implemented, numerical tests
performed, and results compared. Results are also compared to exact
solutions when available.

In sections 2 a short overview on the Adomian decomposition method and the
Euler nonlinear and linearized equations are given. In sections 3 details of
the proposed ABS and ABS-DG schemes are provided. In section 4 a connection
between the ABS and Runge-Kutta (RK) methods is established for the linear
case. Tests are performed and reported in section 5, while conclusions are drawn
in section 6.

\label{Intro}




\section{Review} In this section a short review on the Adomian decomposition
method and nonlinear and linearized Euler equations are given.
\subsection{The Adomian Decomposition Method }
\label{Adomian}
In the following, a short description of the Adomian decomposition method is
given. For more details we refer to \cite{adm1, adm2, adm3, adm5}. The first
step of the method consists in identifying the differential equations in the
following form,
\begin{equation}
    L(u)+R(u)+N(u) = 0
\label{eq1}
\end{equation}
Where $L$ and $R$ are the linear part of the differential operator, with $L$
being the part that is easily invertible. $N$ is the nonlinear part.

Note that it is not necessary to distinguish the non invertible linear part
from the nonlinear invertible one. Both can be represented by the sum $N+R$ by
a single operator, $N$.

Adomian algorithm considers the solution $u$ as a summation over a series,
\begin{equation}
u=\sum_{n=0}^{\infty} u_{n}.
\label{eq2}
\end{equation}
While operator $N$ is given by the expansion of the \textit{Adomian
polynomials}, $N_n$.
\begin{equation}
N=\sum_{n=0}^{\infty} N_{n}.
\label{eq3}
\end{equation}
where the Adomian polynomial coefficients $N_{n}$ are given by
\begin{equation}
N_{n}=\frac{1}{n!}\frac{\partial^{n}}{\partial\lambda^{n}}\left[
N\left(\sum_{k=0}^{n}\lambda^{k}u_{k} \right) \right]_{\lambda=0}.
\label{eq4}
\end{equation}
Finally, $u_{n}$ terms are computed as
\begin{equation}
    u_{n+1}=L^{-1}(N_{n})
    \label{eq4}
\end{equation}
\subsection{The non-conservative Euler equations}
\label{Non-Linear}
The two-dimensional compressible inviscid flow equations are given by the
Euler equations. Formulated relative to a Cartesian $(x, y)$
coordinate system, and in a primitive variable form as
\begin{equation}
\left\{ \begin{array}{lr}
    \frac{\partial \rho}{\partial t} +\partial_{x} (\rho u) +\partial_{y} (\rho v) = 0  \\
\frac{\partial u}{\partial t} +u\partial_{x} (u) +v\partial_{y} (u) + \frac{1}{\rho}\partial_{x}p = 0 \\
\frac{\partial v}{\partial t} +u\partial_{x} (v) +v\partial_{y} (v) + \frac{1}{\rho}\partial_{y}p = 0 \\
\frac{\partial p}{\partial t} +u\partial_{x} (p) +v\partial_{y} (p) + \gamma p(\partial_{x} (u) +\partial_{y} (v)) = 0 \\
        \end{array} \right.
\label{nonlinearEuler}
\end{equation}
where $\rho$ and $p$ denote the averaged density and pressure of the fluid
respectively and $E = e+\frac{1}{2}( u_{\alpha}u_{\alpha} )$ is the total
energy per unit mass, with $e$ being the internal energy per unit mass.
$u_{\alpha}$ is the averaged velocity of the fluid in direction $x_{\alpha}$
and $\delta_{\alpha i}$ is the Kronecker delta.

The equation set is closed with the addition of the ideal gas state equation.
i.e. $p = \rho (\gamma -1) e$, where $\rho$ is the density and $\gamma$ is
the ratio of specific heats.  Solutions to the resulting set of equations
are defined on a fixed spatial computational domain $\Omega$.

Euler equations can also be written in vector or matrix form,
\begin{equation}
    \label{eq5}
    \frac{\partial Q }{\partial t} + \mathcal{A}\frac{\partial Q }{\partial x}
    + \mathcal{B}\frac{\partial Q' }{\partial y} = S.
\end{equation}
Where, $Q(x,y,t) = (\rho, u_1, u_2, p)^t$ and
\begin{align*}
& \mathcal{A}(x,y,t) = \begin{pmatrix}
                u & \rho & 0 & 0 \\
                0 & u & 0 & 1/\rho \\
                0 & 0 & u & 0 \\
                0 & \rho c^2 & 0 & u \\
                 \end{pmatrix} &
\mathcal{B}(x,y,t) = \begin{pmatrix}
                v & 0 & \rho & 0 \\
                0 & v & 0 & 0 \\
                0 & 0 & v & 1/\rho \\
                0 & 0 & \rho c^2 & v \\
                 \end{pmatrix}. &
\end{align*}

\subsection{The Linearized Euler Equations}
\label{Governing Equations}

LEE system is obtained after a linearization around a mean flow. This is
achieved by assuming solution is composed of a mean and a perturbation part:
$Q(x,y,t) = Q_0(x,y,t) + Q'(x,y,t)$ (and also for the source term). Additionally, the mean flow values are assumed to satisfy
\begin{align}
    \frac{\partial Q_0 }{\partial t} + \mathcal{A}_0\frac{\partial Q_0 }{\partial x} +
    \mathcal{B}_0\frac{\partial Q_0 }{\partial y} = S_0.
\end{align}
Inserting previous equations into the Euler system, and performing
nondimensionalization, it yields \cite{Blom2003_thesis},
\begin{equation}
    \label{eq5}
    \frac{\partial Q' }{\partial t} + \mathcal{A}_0\frac{\partial Q' }{\partial x} +
    \mathcal{B}_0\frac{\partial Q' }{\partial y} + \mathcal{A}' \dfrac{\partial
    Q_0}{\partial x} + \mathcal{B}' \dfrac{\partial Q_0}{\partial y} = S'.
\end{equation}
Where $Q'=\left( \rho', u', v', p' \right)^t$ and
\begin{align*}
\label{eq6}
& \mathcal{A}_0(x,y,t) = \begin{pmatrix}
                M_1 & 1 & 0 & 0 \\
                0 & M_1 & 0 & 1 \\
                0 & 0 & M_1 & 0 \\
                0 & 1 & 0 & M_1 \\
                 \end{pmatrix} &
\mathcal{B}_0(x,y,t) = \begin{pmatrix}
                M_2 & 0 & 1 & 0 \\
                0 & M_2 & 0 & 0 \\
                0 & 0 & M_2 & 1 \\
                0 & 0 & 1 & M_2 \\
                 \end{pmatrix}. &
\end{align*}
Subscript $0$ in matrices $\mathcal{A}_0, \mathcal{B}_0$, denotes exclusive dependency to mean flow
values. Whereas primes in  $\mathcal{A}'$ and $\mathcal{B}'$ denote evaluation with perturbed
variables.  Those last perturbation matrices can be considered negligible if
mean flow spatial derivatives are moderated,
\begin{align}
   & \mathcal{A}' \dfrac{\partial Q}{\partial x} = 0  & \mathcal{B}' \dfrac{\partial Q}{\partial y} = 0 &
\end{align}
Finally, no source terms will be assumed $ S = S' = 0 $.

In this paper time independent and constant in space  mean values are
considered for simplicity.  Therefore, $\mathcal{A}_0$ and $\mathcal{B}_0$ matrices are
constant. Nonetheless, one could consider space dependent matrices without lost
of generality.

\vspace{2 cm} In summary, the following governing equations are yielded
\begin{equation}
    \label{eq5Simplified}
    \frac{\partial Q'}{\partial t} +
    \mathcal{A}_0\frac{\partial Q'}{\partial x} +  \mathcal{B}_0\frac{\partial Q'}{\partial
    y} = 0
\end{equation}

For convenience we drop
prime symbol from $Q'$.

\section{Adomian Based Schemes (ABS) }

In this paper we propose a cost-effective numerical scheme to solve LEE (eq.
\ref{eq5Simplified}). Method is based on the Adomian decomposition technique
for time discretization, and DG techniques (which can be replaced by other
techniques) for the space discretization. The scheme is assessed by comparison
to Runge-Kutta DG method, and exact solutions.  The DG space discretization is
implemented following the approach proposed by Shu \cite{Shu1998}. For other
classical DG discretization options, the reader is pointed to Cockburns paper
\cite{Cockburn1999}. Details are not provided in this paper since it is a
well-known method. However, the application of Adomian decomposition for time
discretization is described in detail, since it is the main contribution of the
paper.

\subsection{The ABS Scheme Derivation}
\label{ABS1}
Even if the proposed ABS scheme is applied and assessed for the LEE, the scheme
is derived for the general case of nonlinear Euler equations. This is
motivated by the fact that some useful properties simplify the scheme
formulation for both LEE and the nonlinear case.

To apply the Adomian decomposition technique described in section
\ref{Adomian} to Euler equations \eqref{nonlinearEuler}, we propose to set
\begin{center}
$L = \frac{\partial}{\partial t}$ ~~~and then~~~ $L^{-1}= \int_{0}^{t}$
\end{center}
\noindent
Referring by $F_{x}$ and $F_{y}$ to the $x$ and $y$ space derivatives, operator
$N$ is given by
\begin{center}
$N(Q) = F_{x}(Q) + F_{y}(Q)$
\end{center}
with $Q$ being the primitive variables
\begin{align}
Q(x,y) =
\begin{pmatrix}
    \rho \\
    u \\
    v \\
    p
\end{pmatrix}.
\end{align}
Applying decomposition from \eqref{eq2}, in vector notation,
\begin{align}
    \label{eq:Qdecomposition}
    Q = \sum_{k=0}^{\infty}Q_k,
\end{align}
The Adomian coefficients from \eqref{eq3} can be written as
\begin{equation}
    N_{n}(x,y,t)=\frac{1}{n!}\frac{\partial^{n}}{\partial\lambda^{n}}\left[
    F_{x} \left(\sum_{k=0}^{n}\lambda^{k}Q_{k} \right) +F_{y}
    \left(\sum_{k=0}^{n}\lambda^{k}Q_{k} \right) \right]_{\lambda=0}.
\label{Adomnolinearterm}
\end{equation}
Therefore, substituting into the governing equations,
\begin{gather}
    \label{AdomianNonLinearterms}
    N_{n}(x,y,t)=
    = \frac{1}{n!}\frac{\partial^{n}}{\partial\lambda^{n}}\left[
    \begin{array}{l}
        \partial_{x}
        \left((\sum_{k=0}^{n}\lambda^{k}\rho_{k})(\sum_{k=0}^{n}\lambda^{k}u_{k})\right)
        \\[-5 pt]
        \hspace{1 cm}+\partial_{y}
        \left((\sum_{k=0}^{n}\lambda^{k}\rho_{k})(\sum_{k=0}^{n}\lambda^{k}v_{k})\right)
        = 0  \\[10 pt]
        (\sum_{k=0}^{n}\lambda^{k}u_{k})\partial_{x}(\sum_{k=0}^{n}\lambda^{k}u_{k})
        +(\sum_{k=0}^{n}\lambda^{k}v_{k})\partial_{y}(\sum_{k=0}^{n}\lambda^{k}u_{k})
        + \\[-2 pt]
        \hspace{1 cm}
        \left(\dfrac{1}{\sum_{k=0}^{n}\lambda^{k}{\rho}_{k}}\right)
        \partial_{x} \left(\sum_{k=0}^{n}\lambda^{k}p_{k}\right) = 0 \\[10 pt]
        (\sum_{k=0}^{n}\lambda^{k}u_{k})\partial_{x}(\sum_{k=0}^{n} \lambda^{k}v_{k})
        +(\sum_{k=0}^{n}\lambda^{k}v_{k})\partial_{y}(\sum_{k=0}^{n}\lambda^{k}v_{k})
        + \\[-2 pt]
        \hspace{1 cm} \left(\dfrac{1}{\sum_{k=0}^{n}\lambda^{k}{\rho}_{k}}\right)
        \partial_{y} \left(\sum_{k=0}^{n}\lambda^{k}p_{k}\right)=
        0 \\[10 pt]
        (\sum_{k=0}^{n}\lambda^{k}u_{k})\partial_{x}
        (\sum_{k=0}^{n}\lambda^{k}p_{k})
        +(\sum_{k=0}^{n}\lambda^{k}v_{k})\partial_{y}
        (\sum_{k=0}^{n}\lambda^{k}p_{k}) + \\[-5 pt]
        \hspace{1 cm}
        \gamma
        (\sum_{k=0}^{n}\lambda^{k}p_{k})\left(\partial_{x}
        (\sum_{k=0}^{n}\lambda^{k}p_{k}) +\partial_{y}
        (\sum_{k=0}^{n}\lambda^{k}v_{k})\right) = 0
    \end{array} \right]_{\lambda=0}.
\end{gather}
And $Q_{n+1}$ terms are computed recursively by \eqref{eq4},
\begin{equation}
    \label{InvertibleLinearTerm}
    Q_{n+1}(x,y,t)= \int_{0}^{t}N_{n}(x,y,t).
\end{equation}


Let's expand each equation of the vector $N_n = (A_{n}, B_{n}, C_{n},
D_{n})^{t}$, corresponding each component to an equation of the Euler
system. We will first derive the term $A_{n}$ corresponding to the continuity
equation. Then $B_{n}, C_{n}, D_{n}$ will be obtained in a similar way.

From formula \eqref{AdomianNonLinearterms} $A_{n}$ is given by:
\begin{align*}
    & A_{n}(x,y,t)= \frac{1}{n!}\frac{\partial^{n}}{\partial\lambda^{n}}\left[
    \frac{\partial}{\partial x}
    \left((\sum_{k=0}^{n}\lambda^{k}\rho_{k})(\sum_{k=0}^{n}\lambda^{k}u_{k})
    \right) +\frac{\partial}{\partial
    y}\left((\sum_{k=0}^{n}\lambda^{k}\rho_{k})(\sum_{k=0}^{n}\lambda^{k}v_{k})
    \right) \right]_{\lambda=0} &
\end{align*}
Let's develop the first term in the summation. First, the derivative order
is exchanged. Meaning,
\begin{align*}
    & \frac{1}{n!}\frac{\partial^{n}}{\partial\lambda^{n}}\left[
    \frac{\partial}{\partial x} \left( \cdot \right) \right]_{\lambda=0} =
    \frac{1}{n!}\frac{\partial}{\partial x}\left[
    \frac{\partial^{n}}{\partial\lambda^{n}} \left( \cdot \right)
    \right]_{\lambda=0}. &
\end{align*}
Using the \textit{Leibniz} formula, we get
\begin{align*}
    \frac{\partial^{n}}{\partial\lambda^{n}}
    \left(\sum_{k=0}^{n}\lambda^{k}\rho_{k}\sum_{k=0}^{n}\lambda^{k}u_{k}
    \right) & = \sum_{j=0}^{n}
    \binom{j}{k}\frac{\partial^{n-j}}{\partial\lambda^{n-j}}\left(
    \sum_{k=0}^{n}\lambda^{k}\rho_{k}\right)
    \frac{\partial^{j}}{\partial\lambda^{j}}
    \left(\sum_{k=0}^{n}\lambda^{k}u_{k}\right)
\end{align*}
with
\begin{align*}
    & \frac{\partial^{n-j}}{\partial\lambda^{n-j}}
    \left(\sum_{k=0}^{n}\lambda^{k}\rho_{k}\right) |_{\lambda=0} = (n-j)!\,
    \rho_{n-j} &
    \binom{j}{k} = \dfrac{n!}{j! \ (n-j)!} & \\
    & \frac{\partial^{j}}{\partial\lambda^{j}}\left( \sum_{k=0}^{n}
    \lambda^{k} u_{k}\right) |_{\lambda=0} = (j)!\, u_{j} &
\end{align*}
Then
\begin{align}
    \label{eq:convectionExpansion_1}
    \frac{1}{n!}\frac{\partial}{\partial x}\left[
    \frac{\partial^{n}}{\partial\lambda^{n}}
    \left((\sum_{k=0}^{n}\lambda^{k}\rho_{k})(\sum_{k=0}^{n}\lambda^{k}u_{k})
    \right) \right]_{\lambda=0} = \sum_{j=0}^{n}
    \frac{\partial}{\partial x}(\rho_{n-j}u_{j})
\end{align}
Similarly, for the second summation term we have,
\begin{align}
    \label{eq:convectionExpansion_2}
    \frac{1}{n!}\frac{\partial}{\partial y}\left[
    \frac{\partial^{n}}{\partial\lambda^{n}}
    \left((\sum_{k=0}^{n}\lambda^{k}\rho_{k})(\sum_{k=0}^{n}\lambda^{k}v_{k})
    \right) \right]_{\lambda=0} = \sum_{j=0}^{n}\frac{\partial}{\partial
    y}(\rho_{n-j}v_{j}).
\end{align}

Substituting both equations (\ref{eq:convectionExpansion_1} -
\ref{eq:convectionExpansion_2}) into $A_{n}$ equation, we get the final formula
\begin{align}
    A_{n}=-\sum_{j=0}^{n}(\partial_{x}(\rho_{n-j}u_{j}) +
    \partial_{y}(\rho_{n-j}v_{j}))
\end{align}

We obtain a similar formula for $D_{n}$.

For $B_{n}$ and $C_{n}$, first develop
$\frac{1}{\sum_{k=0}^{n}\lambda^{k}\rho_{k}}$ as power series (note that
$\lambda$ can be considered close to zero since we are concerned by the limit)

\begin{center}
$\frac{1}{\sum_{k=0}^{n}\lambda^{k}\rho_{k}} = \sum_{k=0}^{\infty}\lambda^{k}\widehat{\rho}_{k}$
\end{center}
That is,
\begin{align*}
    1 & =(\sum_{k=0}^{n}\lambda^{k}\rho_{k})(
      \sum_{k=0}^{\infty}\lambda^{k}\widehat{\rho}_{k}) \\
     & =\rho_{0}\widehat{\rho}_{0} +
     \sum_{k=1}^{\infty}(\sum_{j=0}^{k}\rho_{j}\widehat{\rho}_{k-j})\lambda^{k}
\end{align*}
Then we obtain the flowing recursive formula
\begin{align}
    \nonumber
    \left\{ \begin{array}{lr}
    \widehat{\rho}_{0}=\frac{-1}{\rho_{0}} \\
    \label{eq:rho_hat_definition}
    \widehat{\rho}_{k}=\frac{-1}{\rho_{0}}\sum_{j=1}^{k}\rho_{j}
    \widehat{\rho}_{k-j},  \qquad \text{for} \ k=1,\ldots, \infty
    \end{array} \right.
\end{align}
such that we can perform a change of variable in system
\eqref{AdomianNonLinearterms} and set,
\begin{center}
    $ \left(\dfrac{1}{\sum_{k=0}^{n}\lambda^{k}{\rho}_{k}}\right) \partial_{x}
     \left(\sum_{k=0}^{n}\lambda^{k}p_{k}\right) = \left(\lambda^k
    \hat{\rho}_k \right) \partial_{x}
    \left(\sum_{k=0}^{n}\lambda^{k}p_{k}\right) $ \\
    $ \left(\dfrac{1}{\sum_{k=0}^{n}\lambda^{k}{\rho}_{k}}\right) \partial_{x}
    \left(\sum_{k=0}^{n}\lambda^{k}p_{k}\right) = \left(\lambda^k \hat{\rho}_k
    \right) \partial_{x} \left(\sum_{k=0}^{n}\lambda^{k}p_{k}\right). $
\end{center}

And therefore using the same simplifications as for $A_{n}$, we finally
obtain the following formula for $N_{n}$
\begin{equation}
N_{n}(x,y,t) = \left\{ \begin{array}{lr}
     A_{n}=-\sum_{j=0}^{n}(\partial_{x}(\rho_{n-j}u_{j}) + \partial_{y}(\rho_{n-j}v_{j}))  \\
	 B_{n}=-\sum_{j=0}^{n}(u_{n-j}\partial_{x}u_{j} + v_{n-j}\partial_{y}u_{j} + \widehat{\rho}_{n-j}\partial_{x}p_{j}) \\
	 C_{n}=-\sum_{j=0}^{n}(u_{n-j}\partial_{x}v_{j} + v_{n-j}\partial_{y}v_{j} + \widehat{\rho}_{n-j}\partial_{y}p_{j}) \\
	 D_{n}=-\sum_{j=0}^{n}(u_{n-j}\partial_{x}p_{j} + v_{n-j}\partial_{y}p_{j} + \gamma p_{j} (\partial_{x}u_{n-j} + \partial_{y}v_{n-j})) \\
        \end{array} \right.
\label{FormulaForRn}
\end{equation}

In practice, time integration \eqref{InvertibleLinearTerm} is in general not
easy to compute, or at least in a very accurate way. Indeed, Adomian series
coefficients are polynomials in time whose computed coefficients need to be
stored. The following Theorem remedy to this problem and allows a
systematic and exact time integration by a simple multiplication by time $t$.

\textit{\textbf{Theorem}:  In formula (\ref{FormulaForRn}) the
expression of $N_{n}$ can be expressed as
\begin{equation}
    N_{n}(x,y,t) = t^{n}\overline{N}_{n}(x,y)
    \label{RnSimplifiedFormula}
\end{equation}
Where $\overline{N}_{n}(x,y)$, is an expression depending only on $x$ and $y$}
\vskip 0.25cm

\textit{Proof}: We establish the proof by induction. We will do it only for the
first momentum equation terms ($B_{n}$) and the others are obtained in a
similar way.

Note that if equation \eqref{RnSimplifiedFormula} is fulfilled, a similar
relation is hold by field variables,
\begin{equation}
Q_{n+1}(x,y,t) = \int_{0}^{t}N_{n}(x,y,t)dt=
\frac{t^{n+1}}{(n+1)}\overline{N}_{n}(x,y) = t^{n+1}\overline{Q}_{n}(x,y).
\label{eq16}
\end{equation}
Note also that in the recursive formula of $\widehat{\rho}_{k}$ the sum of
the indexes is always equal to $k$, this implies that $\widehat{\rho}_{k}$
satisfies the formula \ref{RnSimplifiedFormula} as long as $\rho_{k}$
satisfies it.
\begin{equation}
    \widehat{\rho}_{k}(x,y,t) = t^{k}\widehat{\overline{\rho}}_{k}(x,y)
    \label{rohat}
\end{equation}

Now to initialize the recursive proof, let's verify the relation  for
$B_{o}$ and $B_{1}$. From \eqref{FormulaForRn} We have
\begin{center}
$B_{0}(x,y,t)=-(u_{0}\partial_{x}u_{0} +v_{0}\partial_{x}u_{0}  +
\frac{1}{\rho_{0}}\partial_{x}p_{0}) =t^{0} \overline{B}_{0}(x,y)$.
\end{center}
Then a new state variable term is computed with
\eqref{FormulaForRn},
\begin{equation}
Q_{1}(x,y,t)= \int_{0}^{t}N_{0}(x,y)=tN_{0}(x,y).
\label{eq12}
\end{equation}
And therefore a new Adomian polynomial,
\begin{center}
$B_{1}=-(u_{1}\partial_{x}u_{0} +v_{1}\partial_{x}u_{0}  + \widehat{\rho}_{1}\partial_{x}p_{0} +
u_{0}\partial_{x}u_{1} +v_{0}\partial_{x}u_{1}  + \widehat{\rho}_{0}\partial_{x}p_{1}) $
\end{center}
Recall from \eqref{rohat},
\begin{center}
    $\widehat{\rho}_{1} = t\widehat{\overline{\rho}}_{1}(x,y)$
\end{center}
Substituting $\widehat{\rho}_{1}$ and $Q_{1}$ into the expression of $B_{1}$,
it is easy to verify that $B_{1}(x,y,t)=t\overline{B}_{1}(x,y)$

Now assume relation (\ref{RnSimplifiedFormula}) is valid till index $n$
and lets proof it for $n+1$. We proceed exactly as for $N_{1}$, since the
relation is valid for order $n$ we have
\begin{center}
    $N_{n}(x,y,t) = t^{n}\overline{N}_{n}(x,y)$
\end{center}
Then,
\begin{equation}
    Q_{n+1}(x,y,t) = \int_{0}^{t}N_{n}(x,y,t) =
    \frac{t^{n+1}}{n+1}\overline{N}_{n}(x,y) = t^{n+1}\overline{Q}_{n}(x,y)
    \label{Qformulat}
\end{equation}

By substituting this expression in equation \eqref{FormulaForRn} we obtain for
the first momentum $B_{n+1}$.
\begin{align*}
    B_{n+1} & = \sum_{j=0}^{n} u_{n+1-j}\partial_{x}u_{j} + v_{n+1-j}\partial_{y}u_{j} +
    \widehat{\rho}_{n+1-j}\partial_{x}p_{j}  \\
    & = \sum_{j=0}^{n} \left[ \left(t^{n+1-j}\overline{u}_{n-j}(x,y)\right)
    \left(\partial_{x} t^{j} \overline{u}_{j-1}(x,y) \right) + \left( t^{n+1-j}
    \overline{v}_{n-j}(x,y) \right) \left( \partial_{y}
    t^{j}\overline{u}_{j-1}(x,y) \right) \right. \\
    & \left. \hspace{1 cm} + \left( t^{n+1-j} \widehat{{\rho}}_{n+1-j}
    \right) \left( \partial_{x} t^{j}\overline{p}_{j-1}(x,y)\right) \right] \\
    & = \sum_{j=0}^{n} t^{n+1}\left[ \overline{u}_{n-j}(x,y)\partial_{x}\overline{u}_{j-1}(x,y)
    + \overline{v}_{n-j}(x,y)\partial_{y}\overline{u}_{j-1}(x,y)+
    \widehat{{\rho}}_{n+1-j}\partial_{x}\overline{p}_{j-1}(x,y)
    \right]
\end{align*}
Which gives,
\begin{equation*}
    B_{n+1}=  -t^{n+1}\sum_{j=0}^{n}
    \left[\overline{u}_{n-j}(x,y)\partial_{x}\overline{u}_{j-1}(x,y) +
    \overline{v}_{n-j}(x,y)\partial_{y}\overline{u}_{j-1}(x,y)+
    \widehat{{\rho}}_{n+1-j}\partial_{x}\overline{p}_{j-1}(x,y)\right]
\end{equation*}
This is a monomial of degree $n+1$ in time with a coefficient depending only
on $x,y$,  therefore $B_{n+1}$ can be written as
\begin{center}
$B_{n+1}(t,x,y)=t^{n+1}\overline{B}_{n}(x,y)$
\end{center}
Which achieves the \textit{proof} of the Theorem.

\qed

\subsection{The ABS formula}
\label{ABS-formula}
Using equation \eqref{Qformulat} we can derive a formula for $Q_{n+1}$ that
doesn't require any time integration,
\begin{align}
    \nonumber
    Q_{n+1}(x,y,t) & = \int_{0}^{t}{N_{n}(t,x,y)dt} = \int_{0}^{t}{t^{n}\overline{N}_{n}(x,y)dt} \\
    \nonumber
    & = \frac{t^{n+1}}{n+1}\overline{N}_{n}(x,y) = \frac{t}{n+1}\left[ t^{n}\overline{N}_{n}(x,y)\right] \\
    & = \frac{t}{n+1}N_{n}(t,x,y)
    \end{align}
Substituting the expression of $N_{n}$ we obtain
\begin{equation}
Q_{n+1}(x,y,t)=\left\{ \begin{array}{lr}
     \frac{-t}{n+1}\left[ \sum_{j=0}^{n}(\partial_{x}(\rho_{n-j}u_{j}) + \partial_{y}(\rho_{n-j}v_{j}))\right]   \\
	 \frac{-t}{n+1}\left[\sum_{j=0}^{n}(u_{n-j}\partial_{x}u_{j} + v_{n-j}\partial_{y}u_{j} + \widehat{\rho}_{n-j}\partial_{x}p_{j})\right] \\
	 \frac{-t}{n+1}\left[\sum_{j=0}^{n}(u_{n-j}\partial_{x}v_{j} + v_{n-j}\partial_{y}v_{j} + \widehat{\rho}_{n-j}\partial_{y}p_{j})\right] \\
	\frac{-t}{n+1}\left[ \sum_{j=0}^{n}(u_{n-j}\partial_{x}p_{j} + v_{n-j}\partial_{y}p_{j} + \gamma (\partial_{x}u_{n-j} + \partial_{y}v_{n-j})p_{j})\right] \\
        \end{array} \right.
\label{Qrecurssiveformula}
\end{equation}

As it can be seen, once $\{Q_{1},\ldots,Q_{n}\}$ are calculated $Q_{n+1}$ is
obtained without any time integration.

Considering space derivatives, the recursive formulation can lead to a fully
analytical formulation, only depending on the initial condition. In literature
some articles can be found where Mathematica is used to directly obtain an
analytical solution to a simple initial condition (see for example \cite{adm5,
Momani2006}). Nevertheless, using symbolic software to obtain a real problem
solution can be very expensive in terms of computational cost, assuming it is
possible to find an analytical expression for the initial condition problem.
Moreover, calculating the derivatives in strong sense requires the initial
condition to be smooth enough.

We hereby propose considering space derivatives in the weak sense and estimate
them numerically. Although any numerical method can be used, in this work space
derivatives are estimated using discontinuous Galerkin (DG) method as proposed
by Shu \cite{Shu1998}. However, different DG methods can be found in, for
instance, Cockburn
\cite{Cockburn1999}. We refer to the obtained (fully discretized) scheme by ABS-DG.

Now and as mentioned in the introduction, detailed formulation of the ABS-DG
scheme will be given for LEE, as it is our target application to assess the
ABS-DG.

\subsection{The ABS-DG for LEE}
For the LEE case, formula \eqref{Qrecurssiveformula} is simplified to
\begin{equation*}
    \label{ABSforLEE}
    Q_{n+1}(x,y,t) = \frac{-t}{n+1} \left[
    \mathcal{A}_0(x,y)\frac{\partial}{\partial x}Q_{n}(x,y,t) +
    \mathcal{B}_0(x,y)\frac{\partial}{\partial y}Q_{n}(x,y,t)\right]
\end{equation*}
\noindent

\subsection{Space discretization}
Space discretization is achieved by applying the DG method for each term
of the ABS series. The procedure is given in \cite{Shu1998}. However, details for
order zero DG is given here, which corresponds to a finite volume scheme, since
results are used in next step for stability analysis.

In that case each $Q_{n}$ term is approximated at the cell center by
\begin{align*}
    Q_{n}(x_{i},y_{i},t) \simeq \frac{1}{\mid S_{i} \mid}\int_{S_{i}}
    Q_{n}(x,y,t) \, \text{d}S,
\end{align*}
where $S_{i}$ is a given cell surface (for a two-dimensional domain). This
leads to
\begin{align*}
    & Q_{n+1}(x_{i},y_{i},t) = \frac{-t}{n+1} \frac{1}{\mid S_{i} \mid}
    \int_{\partial S_{i}} \Big[ {\mathcal{A}_0}\ Q_{n}(x,y,t) \ \eta_{x} & \\
    & \hspace{5.6 cm} + {\mathcal{B}_0}\ Q_{n}(x,y,t)\ \eta_{y} \Big] \
    \text{d} S  &
    \\
    \label{ABS-DG-ZeroOrder}
    & Q^{N}=\sum_{n=0}^{N}Q_{n} &
\end{align*}

Here $Q^{N}$ is the approximated solution, and $N$ corresponds to the index for
which $\vert Q_{N} \vert$ is smaller than a given tolerance.

Approximating fluxes as in the classical finite volume, we obtain the
following ABS-DG zero order numerical scheme,
\begin{equation*}
    \label{eq21}
(ABS-DG) \left\{    
         Q_{0}(x_{i},y_{i})=Q(x_{i},y_{i},0) ~~ \\[-0 pt]
         Q_{n+1}(x_{i},y_{i},t) = \frac{-t}{n+1} \frac{1}{\mid S_{i} \mid} \sum_{j}
        F_{n}^{i,j}  ~~ \\[-10 pt]
         Q^{N}=\sum_{n=0}^{N}Q_{n} ~~
    \right. .
\end{equation*}
The numerical flux could be the Lax-Friedrichs flux approximation,
\begin{gather}
    \begin{align*}
        F^{i,j}_{n} = \dfrac{1}{2}\Big( \overline{\mathcal{A}} \left( Q_{n}^{i} +
        Q_{n}^j \right) \eta_{x}\ + \overline{\mathcal{B}} \left( Q_{n}^{i}
        + Q_{n}^j \right) \eta_{y} \Big) - \dfrac{1}{2} \alpha \Big( \left( Q_{n}^{j} -
        Q_{n}^{i} \right) \eta_{x} + \left( Q_{n}^{j} - Q_{n}^{i} \right) \eta_{y}
        \Big)
    \end{align*} \\[5 pt]
    \begin{align*}
        & \overline{\mathcal{A}} = \dfrac{\mathcal{A}^{i} + \mathcal{A}^{j}}{2} &
        \overline{\mathcal{B}} = \dfrac{\mathcal{B}^{i} + \mathcal{B}^{j}}{2} &
    \end{align*}
\end{gather}

or any other suitable one.

\subsection{Stability Analysis}

We study stability for the one-dimensional linear wave propagation equation,
\begin{align}
    \dfrac{\partial u}{\partial t} + a \dfrac{\partial u}{\partial x} = 0,
\end{align}
using a zero order spatial  ABS-DG formulation. For the rest of this section,
the spatial discretization index $i$ is written as a superindex, and Adomian
iteration $n$ as a subindex. This non-standard notation is intended to
differentiate Adomian iterations from classical finite difference time levels,
where $n$ usually represents current time level and $n+1$ is the time level
after a time increment. Adomian iteration terms ($u_n$) are given by the
following discretization,
\begin{equation}
  \begin{array}{lr}
     u_{0}^i=u(x_{i},0)   \\
     u_{n+1}^{i}(t) = \dfrac{-t}{n+1} \dfrac{1}{2\,h}
        \Big(a\big(u_{n}^{i+1} - u_{n}^{i-1}\big)-\alpha \,
        \big(u_{n}^{i+1} - 2u_{n}^{i} +u_{n}^{i-1}\big) \Big) .
   \end{array}
\label{onedimequation}
\end{equation}
Where the final solution is obtained by the sum of all the $u_{n}$
contributions, namely; $u=\sum_{n}u_{n}$. In fact, round-off errors for each
term ($\varepsilon_n$) are also governed by the same equation. After applying
a Fourier decomposition for each error term, $\varepsilon_{n}^{i}(t) = \beta_{n}
(t)e^{JK_{n}x_{i}}$, the error modes evolution is obtained,
\begin{equation}
  \begin{array}{lr}
     \varepsilon_{n+1}^{i} = \dfrac{-t}{n+1} \dfrac{\beta_n}{2\,h}
        \Big(a\big(e_{n}^{i+1}-e_{n}^{i-1}\big)-
        \alpha \, \big(e_{n}^{i+1} -2u_{n}^{i} +e_{n}^{i-1}\big) \Big).
  \end{array}
  \label{Errorequation}
\end{equation}
From standard stability procedures, the solution is stable respect to
round-off errors as long as their growth is kept bounded. \textit{Von Neumann}
stability approach is used for this study, with a growth rate amplification
$G_{n} = \frac{e_{n +1}^{i}(t)}{e_{n}^{i}(t)}$. When equation
\eqref{Errorequation} is plugged in, the following equation is obtained,
\begin{align}
    \nonumber
    G_{n} &  = \dfrac{-t}{n+1} \, \dfrac{1}{2\,h} \
    \dfrac{1}{\beta_{n}e^{JK_{n}(x_{i})}} \bigg(
    a\Big(\beta_{n}e^{JK_{n}(x_{i}+h)} -
    \beta_{n}e^{JK_{n}(x_{i}-h)}\Big) \\
    \nonumber
    & \hspace{4 cm} - \alpha \Big(\beta_{n}e^{JK_{n}(x_{i}+h)} -
    2\beta_{n}e^{JK_{n}(x_{i})}+\beta_{n}e^{JK_{n}(x_{i}-h)}\Big) \bigg) \\
    \label{eq:Amplification}
    &= \dfrac{-t}{n+1} \, \dfrac{1}{2\,h} \  \bigg( a\Big(e^{JK_{n}(h)} -
    e^{JK_{n}(-h)}\Big) - \alpha \Big(e^{JK_{n}(h)} + e^{JK_{n}(-h)} - 2
    \Big) \bigg) .
\end{align}
But with $\alpha=\dfrac{1}{2}\dfrac{h}{t}$,
\begin{align*}
    & G_{n} = \dfrac{-1}{2(n+1)} \, \ \bigg( \dfrac{t}{h} a\left(e^{JK_{n}(h)}
    - e^{JK_{n}(-h)}\right) - \dfrac{1}{2} \left(e^{JK_{n}(h)} +
    e^{JK_{n}(-h)} - 2 \right) \bigg) . &
\end{align*}
On the other hand we have
\begin{align*}
  \left(e^{\frac{JK_{n}(h)}{2}} - e^{-\frac{JK_{n}(h)}{2}}\right)^{2}
  & =- 4\sin\left(\frac{JK_{n}(h)}{2}\right)^{2} \\
  & = e^{JK_{n}(h)} + e^{JK_{n}(-h)} - 2 \\
  & = 2j\, \sin(JK_{n}(h)) -2
\end{align*}
By substitution in expression \eqref{eq:Amplification} for $G_{n}$, and
setting $r=|\frac{at}{h}|$ and $\theta_{n}=K_{n}h$
\begin{align*}
    G_{n} &= \dfrac{-1}{2(n+1)} \left( 2rj \, \sin\left(\theta\right) + 2 \,
    \sin\left(\dfrac{\theta}{2}\right)^{2} \right) \\
    & = \dfrac{-1}{2(n+1)} \left( 4rj \, \sin\left(\dfrac{\theta}{2}\right) \,
    \cos\left(\dfrac{\theta}{2}\right) + 2 \,
    \sin\left(\dfrac{\theta}{2}\right)^{2} \right) \\
    & \dfrac{-2 \, \sin\left(\dfrac{\theta}{2}\right)}{n+1} \left( rj \,
    \cos\left(\dfrac{\theta}{2}\right) +
    0.5 \, \sin\left(\dfrac{\theta}{2}\right) \right).
\end{align*}
Then
\begin{align*}
    |G_{n}|= \dfrac{2 \, \left|\sin\left(\dfrac{\theta}{2}\right)\right|}{n+1} \left(
    \sqrt{r^{2} \, \cos\left(\dfrac{\theta}{2}\right)^{2} + 0.25 \,
    \sin\left(\dfrac{\theta}{2}\right)^{2}} \right)
\end{align*}

\begin{figure}[!htbp]
\centering
    \includegraphics[scale=0.5]{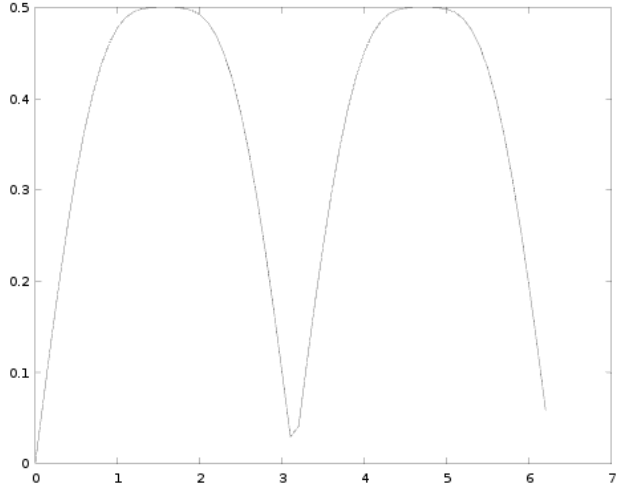} \\
\vskip-0.25cm
\caption{$H(\theta)$ function numerical evaluation}
\label{stabilitycurve}
\end{figure}

Now assume that $r^{2}<\dfrac{(n+1)^{2}}{2}$ meaning that $|r|<\dfrac{n+1}{\sqrt{2}}$, we have
\begin{align*}
    |G_{n}| & < 2 \left|\sin\left(\dfrac{\theta}{2}\right)\right|
    \sqrt{0.5 \, \cos\left(\dfrac{\theta}{2}\right)^{2} +
    \dfrac{0.25}{(n+1)^{2}} \, \sin\left(\dfrac{\theta}{2}\right)^{2}} \\
    & \leq 2 \left| \sin\left(\dfrac{\theta}{2}\right) \right| \sqrt{0.5 \,
    \cos\left(\dfrac{\theta}{2}\right)^{2} + 0.25 \, \sin\left(\dfrac{\theta}{2}\right)^{2}}
\end{align*}
In figure \ref{stabilitycurve} the curve of
$H\left(\overline{\theta}\right)=|\sin\left(\bar{\theta}\right)| \left( \sqrt{0.5 \,
\cos\left(\overline{\theta}\right)^{2} + 0.25 \, \sin\left(\bar{\theta}\right)}
\right)$ as a function of $\overline{\theta}$ is depicted for $0\leq
\overline{\theta}\leq 2\pi$.

We deduce that
\begin{align*}
   r = \left|\dfrac{at}{h}\right| < \dfrac{n+1}{\sqrt{2}}
   ~~~~\text{implies}~~~~ |G_{n}|< 1.
\end{align*}
We conclude that a classical CFL condition is necessary to stabilize the first
Adomian term $(u_{1})$ in the decomposition series \eqref{eq:Qdecomposition}.
For the next terms, as $n$ grows, the condition becomes less restrictive. Note
that the first term requires $1/\sqrt{2}$ instead of the classical $1/2$ for
finite volume with Lax-Friedrichs fluxes, which implies a slight improvement
in stability.

\section{Connections of ABS and RK schemes in the Linear Case}
In this section we will establish a connections between the proposed ABS and the
Runge-Kutta (RK) schemes for the linear case. To solve the ODE
\begin{center}
    $X'=f'(t,x)$
\end{center}
the general form of RK scheme is given by
\begin{center}
    $X_{n+1}=X_{n}+h\sum _{i=1}^{n}c_{i}k_{i}$
\end{center}
where
\begin{align*}
    & k_{1}=f(t_{n},X_{n})  &\\
    & k_{2}=f(t_{n} + \alpha_{2} h,X_{n} + h \beta_{21}k_{1}(t_{n},X_{n}))  &\\
    & k_{3}=f(t_{n} + \alpha_{3} h,X_{n} + h( \beta_{31}k_{1}(t_{n},X_{n})) +
    \beta_{32}k_{2}(t_{n},X_{n})) )  &\\
    \vdots \\
    & k_{m}=f(t_{m} + \alpha_{m} h,X_{n} + h(\sum _{j=1}^{m-1}
    \beta_{mj}k_{j}) &
\end{align*}

Setting $c_{i}=1/i$ and  $\alpha_{j}=0$ for all $j=\{2,\ldots, m\}$ and
$\beta_{kj}=0$ for all $j=\{1,\ldots, m-1\}$ and $k=\{2, \ldots,m\}$ in the
general Runge-Kutta formula we get the ABS scheme. Note that this is not true
in the nonlinear case, it can be easily checked for the Burgers' equations for
instance. The ABS for linear problems appears to be an efficient and a
practical way of applying RK thanks to its recursive formula. Moreover, the
order is dynamic and adaptive for each point of the domain and each timestep.
Being dependent on the remainder of the Adomian series. Therefore there is no
need to fix the order as for the classical RK formulation in advance, and a
maximum accuracy with optimal cost is guaranteed.

\section{Numerical validation and assessment}
To assess the performance of the proposed ABS-DG scheme, two tests are
performed. First a wave propagation is considered, where the
simulation is stopped before the wave reaches the boundary. Hence boundary
effects are avoided. With this appropriate condition, grid convergence is
studied, verifying that the expected spatial order is achieved for various
shape functions polynomial orders. In the second test, non-reflecting
and wall boundary conditions are tested in order to verify they can be
aeroacoustic applications can properly be run.

\subsection{Free-boundary conditions case}
The ABS-DG pressure results are compared to an explicit second-order
Runge-Kutta DG scheme (RK-DG). The test case consists of a Gaussian pulse
centered at the origin propagating for a short period of time (3 nondimensional
time units), such that the simulation is stopped before the wave reaches domain
boundary. Grid dimensions are $30 \times 30$ with a cell edge size of
$0.19$.

Exact solution for pressure in equation \eqref{eq5Simplified} is available
in \cite{adm11} (see $B1-B11$ for details)
\begin{equation*}
    p(x,y,t)=\frac{\varepsilon_{1}}{2\alpha_{1}} \int_{0}^{\infty}\left[
    e^{-\xi^{2}/4\alpha_{1}}cos(\xi t)J_{0}(\xi\eta)\xi \right] d\xi
\end{equation*}
Where $\eta = \left[ (x-Mt)^{2} + y^{2} \right]^{1/2} $ and $J_{0}$ is the zero
order Bessel function. $\alpha_{1}=1/2 \ln{\left(2b\right)}$, $b$ being the half-with
of the Gaussian function. For the performed simulations we set  $\alpha_{1}=1$
and $\varepsilon_{1}=10^{-5}$.

Timestep for reference solution RK-DG is $\Delta t = 0.02$ (equivalently,
\mbox{$\text{CFL} = 0.1$}). For ABS-DG on the other side, \mbox{$\Delta t =
0.5$} (an equivalent of \mbox{$\text{CFL} = 2.5$}. In other words, Adomian
algorithm is restarted after an iteration with $t = 0.5$ is computed. Being the
simulation stopped when Adomian expansion terms are smaller than a set
tolerance value of $10^{-8}$.

\paragraph*{\textbf{Results and discussion}}

Obtained results are discussed next in terms of accuracy and
cost-effectiveness. A grid convergence is also performed for ABS-DG, to
ensure that the right order is obtained.

\paragraph*{\textit{Accuracy assessment}}

\begin{figure}[!htbp]
    \vspace{-2 cm}
    \centering
    \begin{subfigure}{0.95\textwidth}
        \centering
        \includegraphics[width=0.95\linewidth]{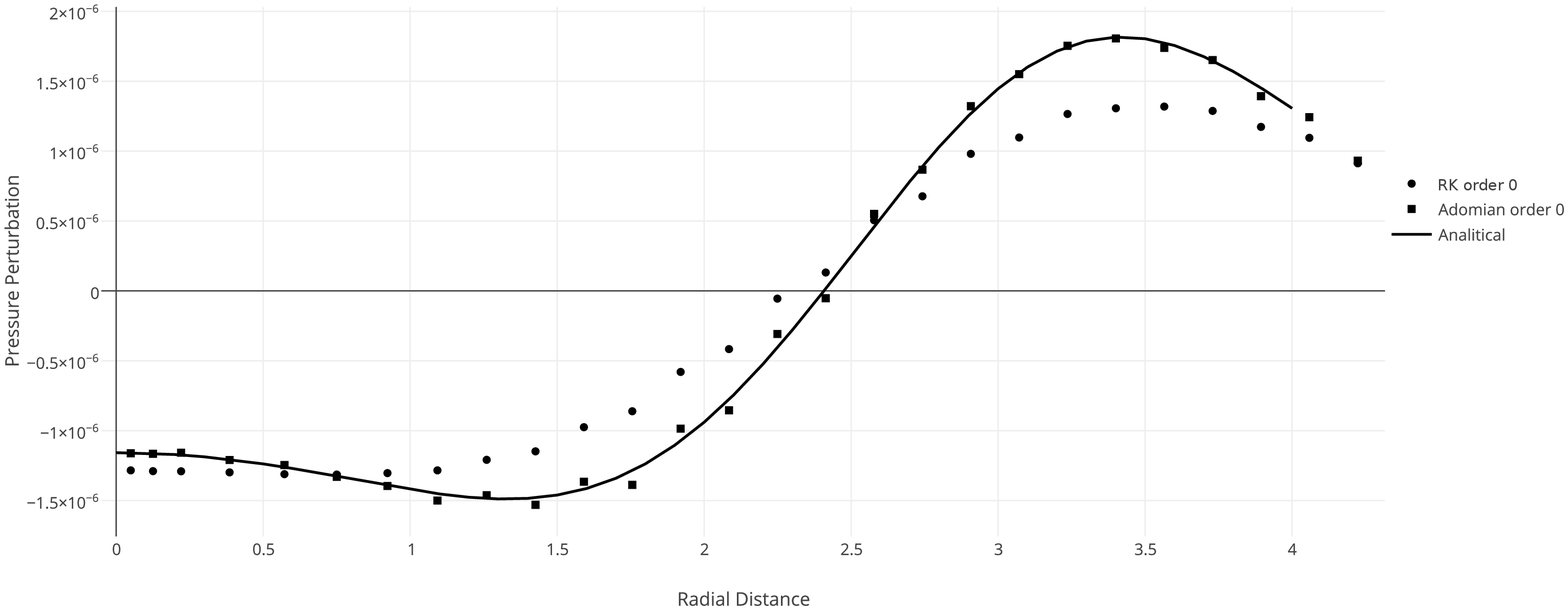}
        \caption{Comparison of the first-order ABS-DG and RK-DG results to the exact solution}
        \label{Order0Comp}
    \end{subfigure}
    \begin{subfigure}[b]{0.95\textwidth}
        \centering
        \includegraphics[width=0.95\linewidth]{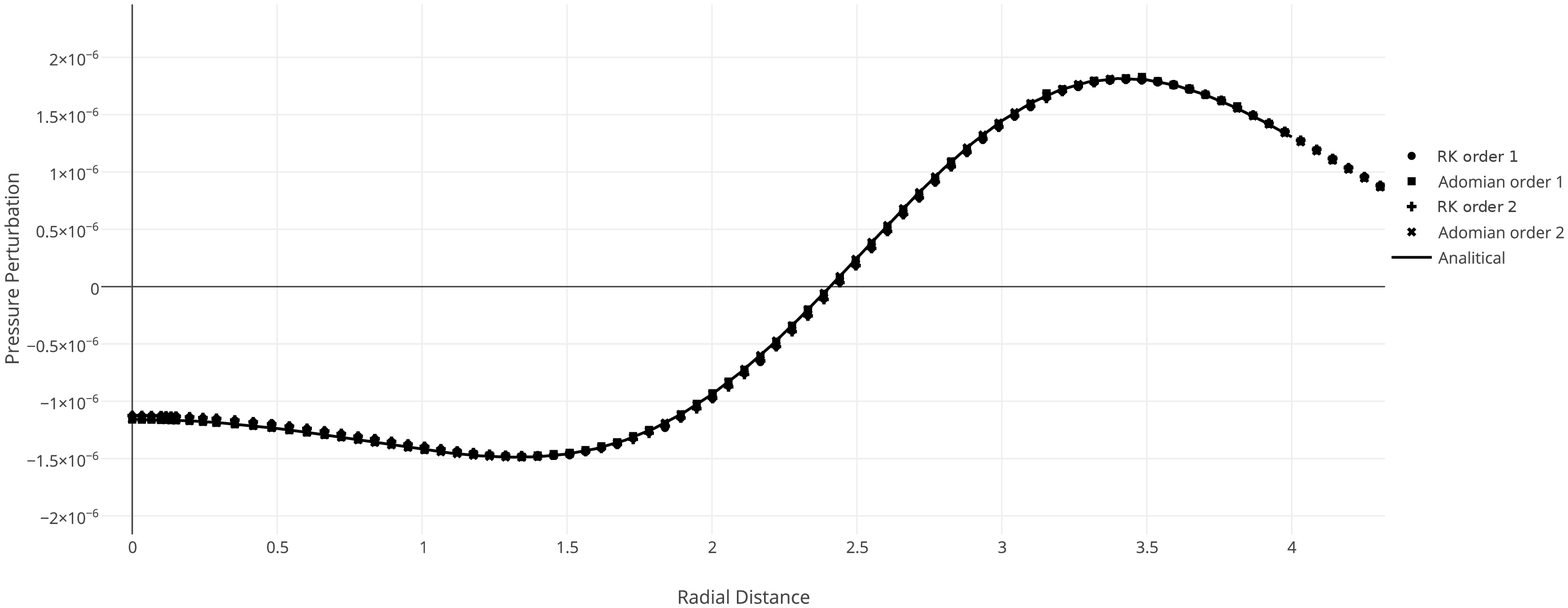}
        \caption{Comparison of the second and third-order ABS-DG and RK-DG
        results to the exact solution } \label{Order1and2}
    \end{subfigure}
    \begin{subfigure}[b]{0.95\textwidth}
    \centering
        \includegraphics[width=0.95\linewidth]{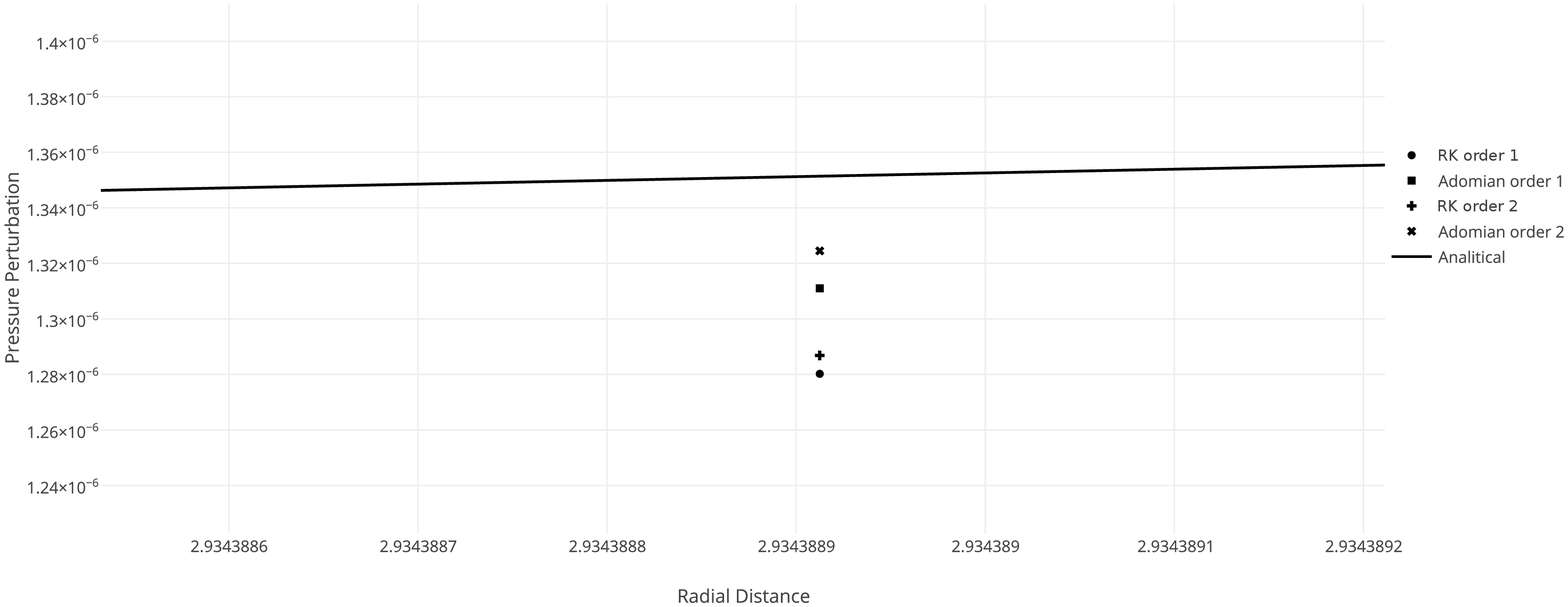}
        \caption{Zoom on ABS-DG, RK-DG and exact solution }
        \label{Order1and2ZOOM}
    \end{subfigure}
    \caption{}
    \label{Adomian-RK-comparison}
\end{figure}

The relative $L^{2}$ error to the exact solution for both RK-DG and ABS-DG
schemes are summarized in table \ref{DGvsABDS} for different spatial orders.
Results show that ABS-DG yields smaller error values than RK-DG. For instance,
ABS-DG first spatial order results are comparable to second-order RK-DG. This
is probably explained by the high accuracy of ABS-DG in time, since the
employed order is dynamic. In other words, the algorithm computes for each cell
the required number of operations such that accuracy satisfies a tolerance at
each point of the domain.  Figure \ref{Adomian-RK-comparison} shows a series
of comparison between ABS-DG and RK-DG respect to the exact solution of the
propagated Gaussian pulse. We can see that ABS-DG results fit better the exact
solution, endorsing results shown in the table.

\begin{table}[ht]
    \caption{Relative $L^{2}$ error for DG vs ABS-DG} 
    \centering 
    \vspace{0.2 cm}
    \begin{tabular}{c c c} 
        \hline\hline 
        Order & DG Method & ABS-DG Method \\ [0.5ex] 
        \hline 
        1 & 2.75E-001 & 6.28E-002 \\ 
        2 & 2.18E-002 & 2.45E-002 \\
        3 & 2.28E-002 & 1.26E-003 \\[1ex] 
        \hline 
    \end{tabular}
    \label{DGvsABDS} 
\end{table}


 \newpage
\paragraph*{\textit{Cost-effectiveness assessment}}

Since both RK-DG and ABS-DG have a similar cost per stage, in order to assess
the cost-effectiveness of the proposed method, the number of computed stage
iterations are compared in table \ref{DGvsABDS_complex} rather than
computational time. In the case of ABS-DG, since different number of iterations
are performed for each cell, the maximum number of Adomian iterations are
counted among all cells. Results indicate that the ABS-DG can reduce the number
of iterations by up to 20 times for the first-order and slightly less
iterations are needed for the third-order.

The total cost of ABS-DG is therefore smaller than RK-DG for the selected
test case. Two arguments in favor of this results are given next. First, it is
well accepted now that in order to obtain high accuracy it is better to
increase the order of the method rather than refining the grid or time spacing
\cite{Kroll2015}. Adomian effectively increases the order in time integration
with each additional iteration. Second, ABS-DG seems to be more stable,
implying less restrictive CFL conditions (recall previous test case was run
with a $\text{CFL = 2.5}$). Despite in this test Adomian requires more stages
per time-iteration, the total cost of the method is given by the product of
the total number of steps and the number of stages per step.

Finally, the important adaptivity property of the ABS-DG scheme is stressed,
which allows for each cell to compute only the required number of iterations to
satisfy a tolerance threshold. Hence, big savings are obtained by avoiding
irrelevant computations on the fly. This property was not reflected neither in
table \ref{DGvsABDS_complex} or in its speedup calculations. As in any adaptive
method, reality is case dependent.

\begin{table}[ht]
    \caption{RK-DG time-iterations VS ABS-DG series terms-iterations for 3
seconds of simulation} 
    \centering 
    \begin{tabular}{c c c} 
        \hline\hline 
        Order & DG(Total time iterations) & ABS-DG(Total series terms iterations) \\ [0.5ex] 
        \hline 
        1 & 300 & 50  \\ 
        2 & 300 & 90  \\
        3 & 300 & 170 \\
        \hline 
    \end{tabular}
    \label{DGvsABDS_complex} 
\end{table}

\paragraph*{\textit{Grid convergence}}

\begin{figure}[ht]
    \centering
        \includegraphics[width=\linewidth]{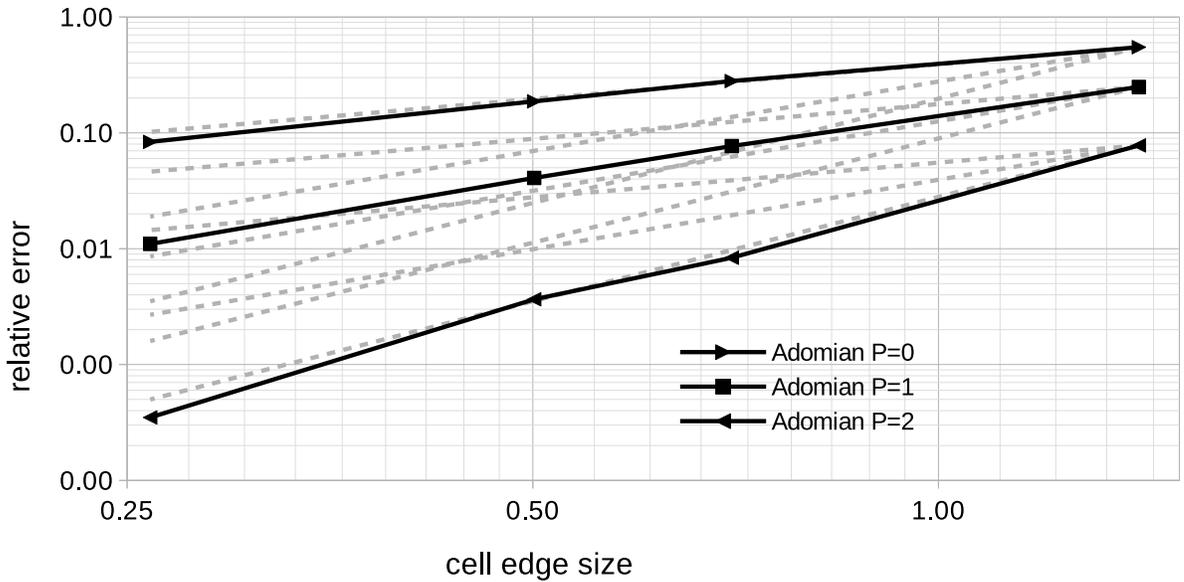}
    \caption{Grid convergence results}
    \label{allOrders}
\end{figure}

\begin{figure}[h]
    \centering
        \begin{subfigure}[t]{0.48\textwidth}
            \centering
            \includegraphics[width=\linewidth]{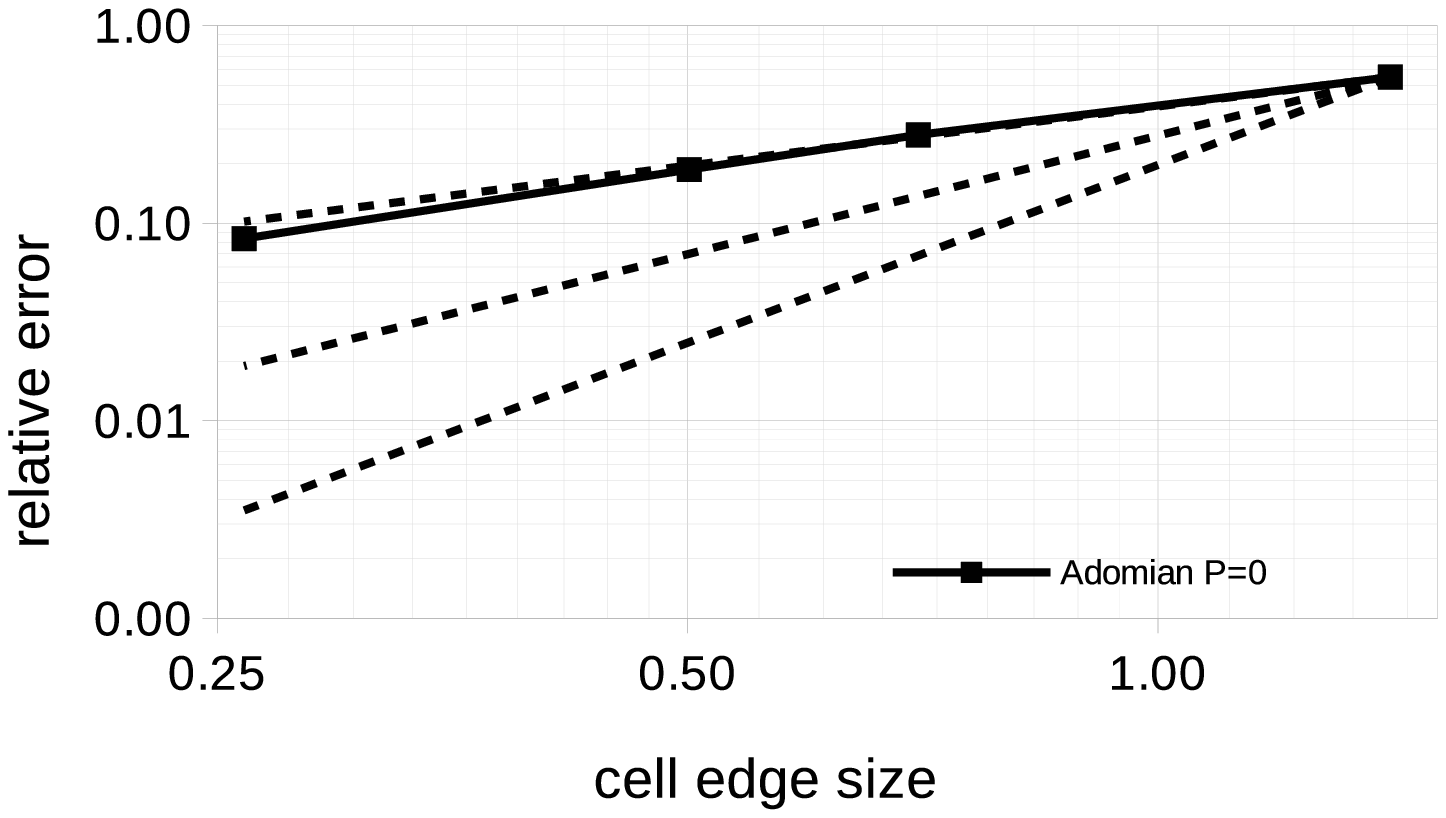}
            \caption{Grid convergence order 0}
            \label{Convergence0}
        \end{subfigure}
        \begin{subfigure}[t]{0.48\textwidth}
            \centering
            \includegraphics[width=\linewidth]{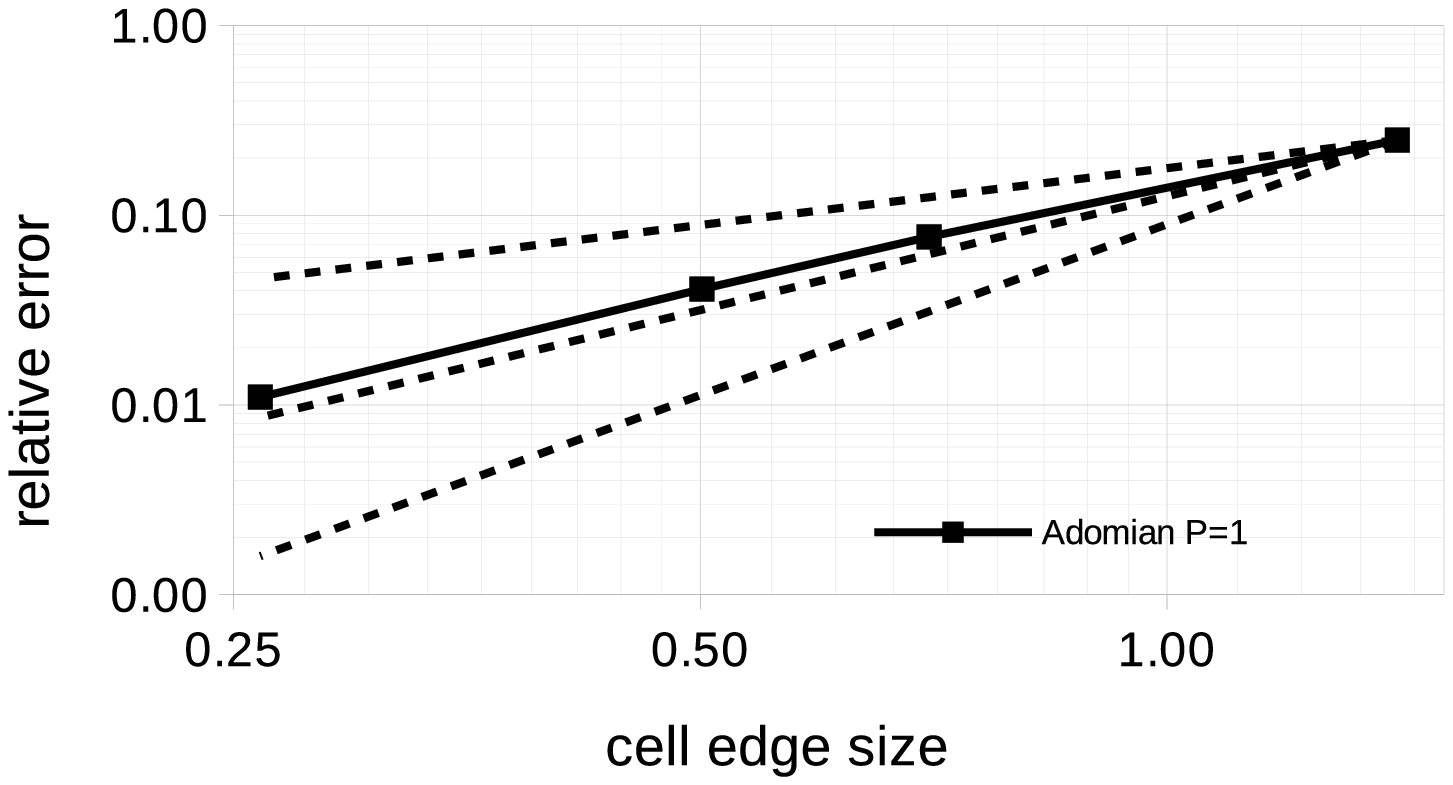}
            \caption{Grid convergence order 1 }
            \label{Convergence1}
        \end{subfigure}
        \begin{subfigure}[t]{0.48\textwidth}
            \includegraphics[width=\linewidth]{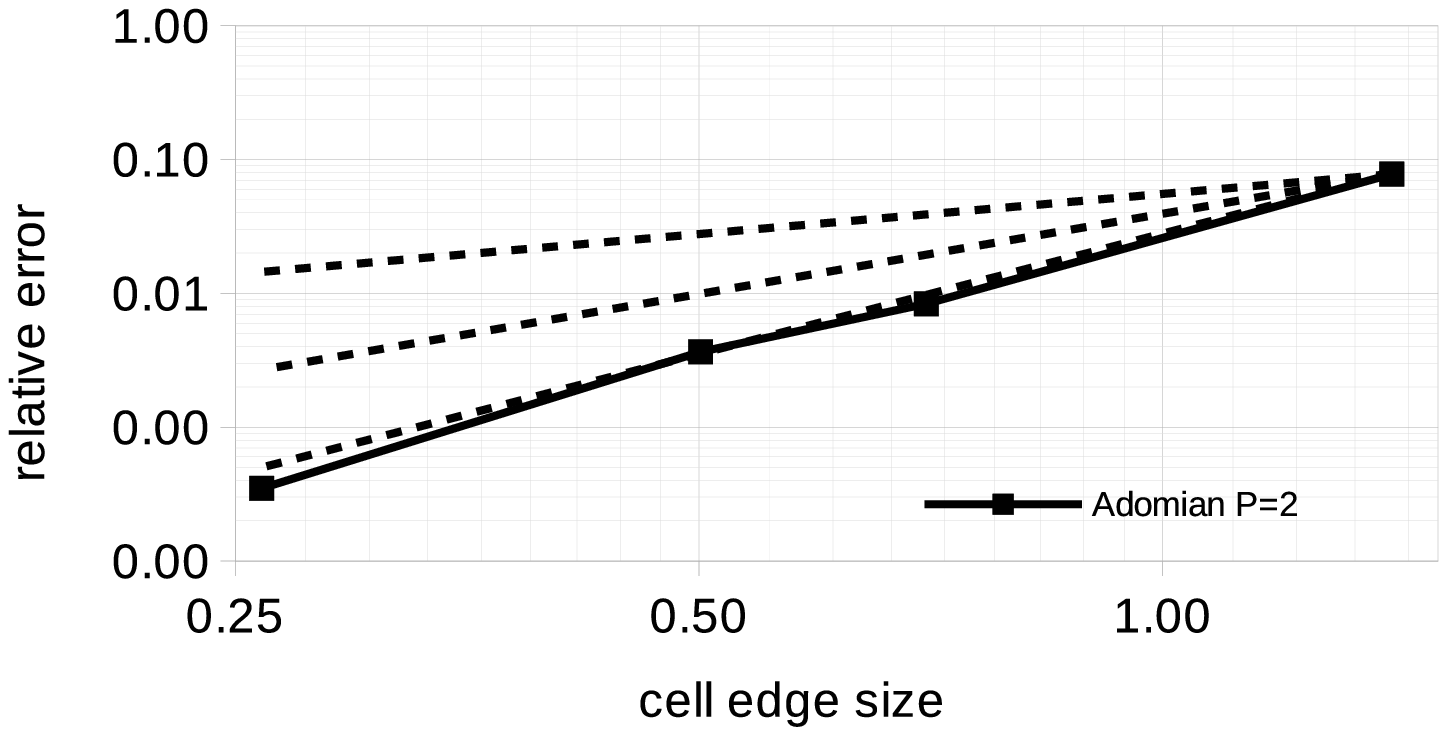}
            \caption{Grid convergence order 2 }
            \label{Convergence2}
        \end{subfigure}
    \caption{}
    \label{fig:ConvergenceAllSeparated}
\end{figure}

To study grid convergence of the ABS-DG method, four different meshes with
different sizes are generated. This is done by selecting different edge size in
SALOME, a tool used to generate meshes on the current work \cite{Bergeaud2010}.
Simulations are stopped at $t=2$, and the relative (to the exact) error is
computed. Table \ref{GridConvergence} shows the errors for different $h$ size
and in figure \ref{allOrders} their logarithmic curves are plotted. For
clarity, each curves is separately shown in figure
\ref{fig:ConvergenceAllSeparated}. In dash line the theoretical order of
convergence is plotted. The solid line represents the numerical solution.
These results demonstrate that we get the right order.

\begin{table}[ht]
    \caption{Relative error magnitude for several Adomian orders at different Mesh sizes} 
    \centering 
    \begin{tabular}{c c c c c} 
        \hline\hline 
        Order & $h=1.4$ & $h=0.7$ & $h=0.5$ & $h=0.26$ \\[0.5ex] 
        \hline 
        1 & 5.49E-001 & 2.80E-001 & 1.87E-001 & 8.36E-002 \\ 
        2 & 2.49E-001 & 7.69E-002 & 4.08E-002 & 1.10E-002 \\
        3 & 7.79E-002 & 8.36E-003 & 3.66E-003 & 3.49E-004 \\[1ex] 
        \hline 
    \end{tabular}
    \label{GridConvergence} 
\end{table}

%
%

%
%

\newpage
\subsection{Tests with boundary conditions}
\begin{figure}[h]
    \centering
    \begin{subfigure}[t]{0.45\textwidth}
        \centering
        \includegraphics[width=\linewidth]{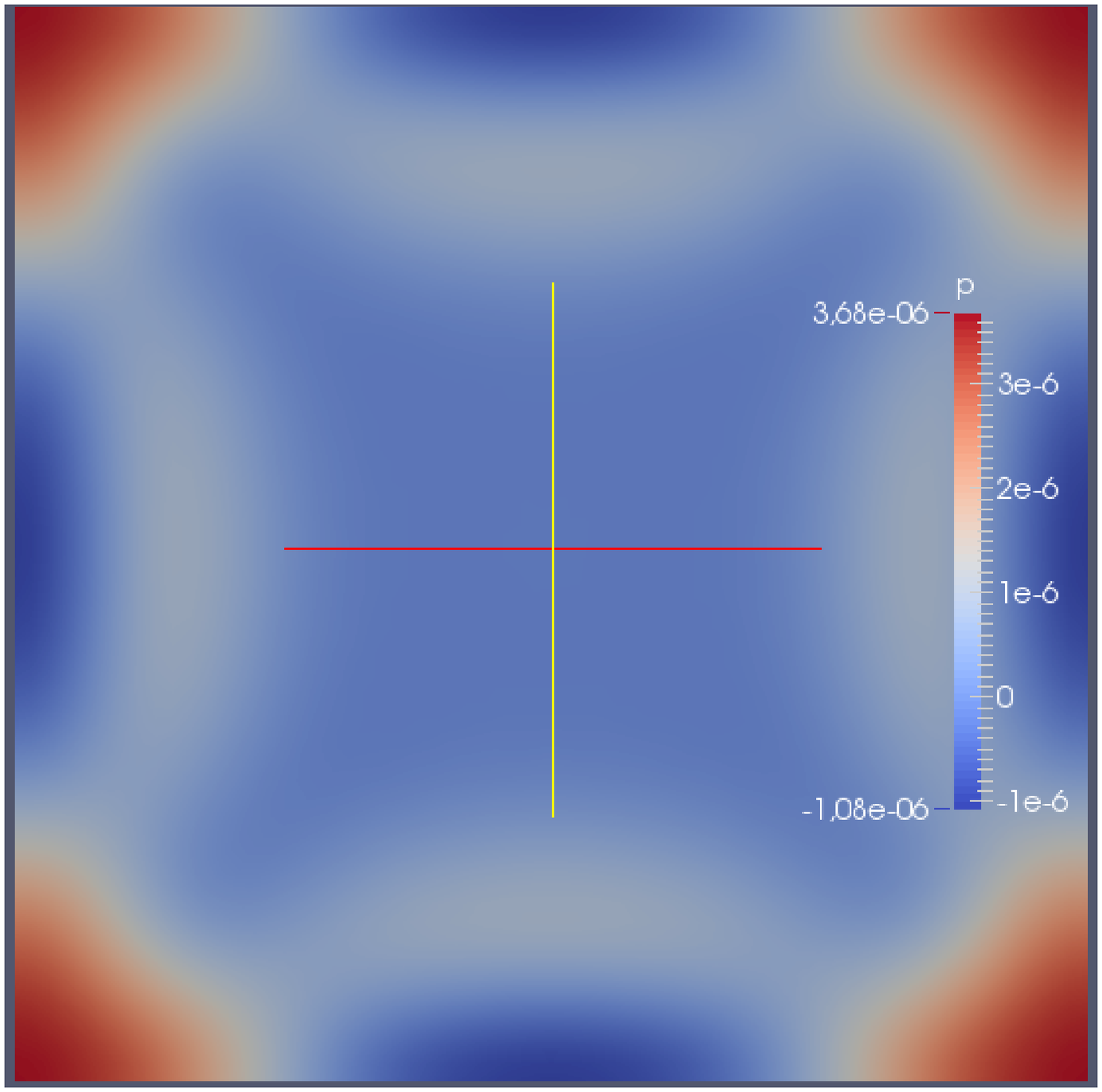}
        \caption{Reference solution}
    \end{subfigure}
    \begin{subfigure}[t]{0.45\textwidth}
        \centering
        \includegraphics[width=\linewidth]{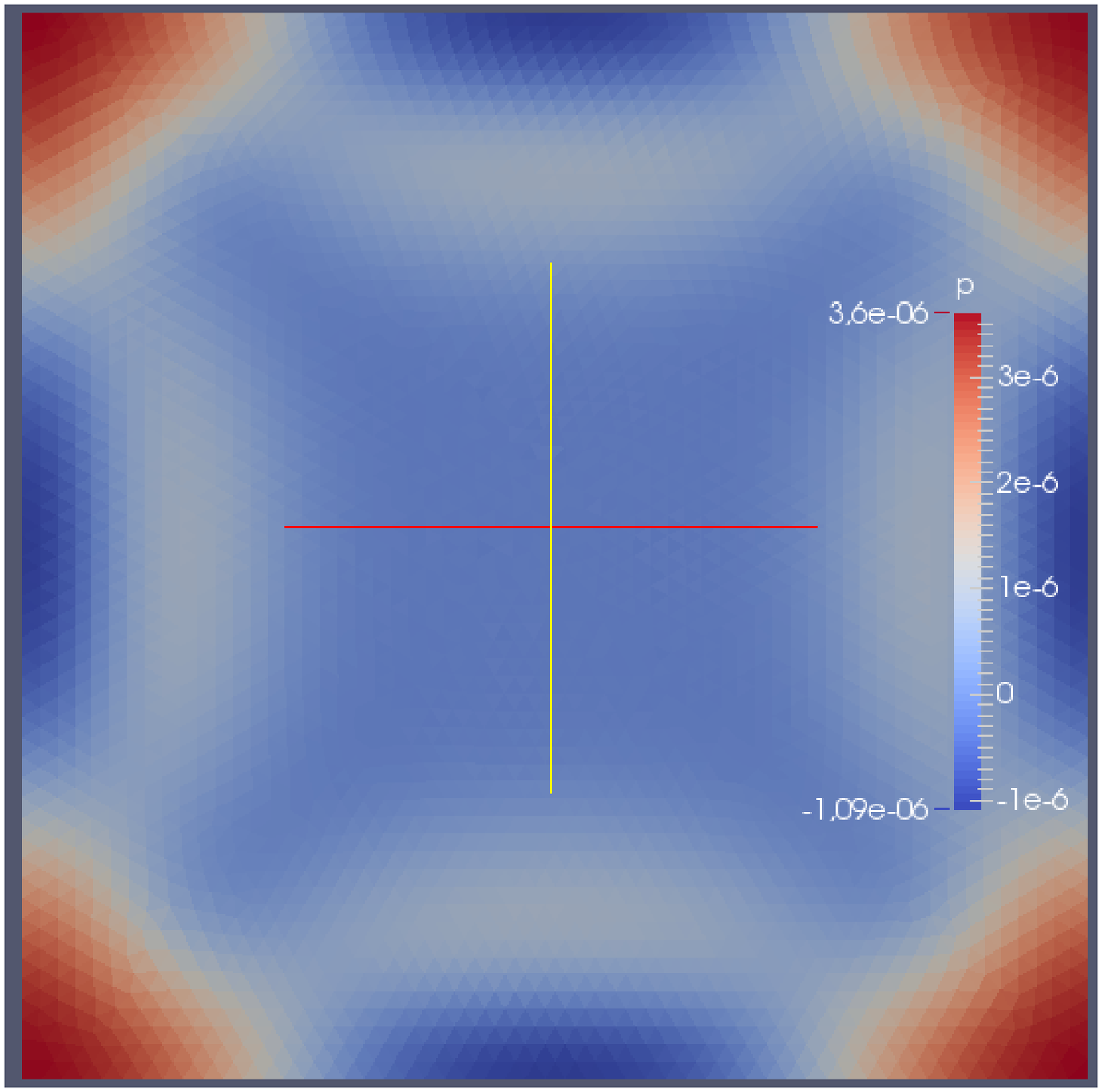}
        \caption{P0 Adomian solution}
    \end{subfigure}
    \begin{subfigure}[t]{0.45\textwidth}
    \centering
        \includegraphics[width=\linewidth]{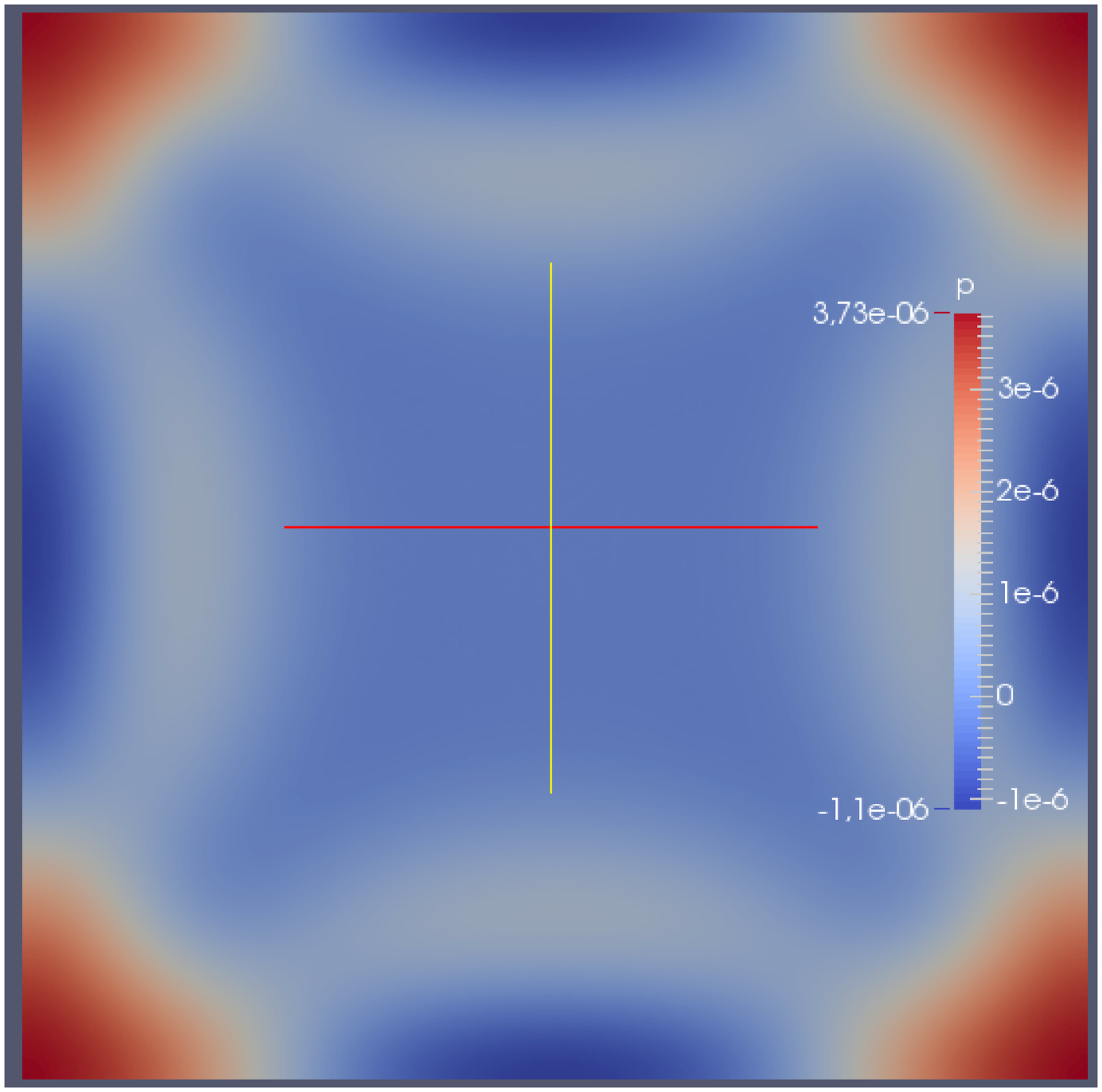}
        \caption{P1 Adomian solution} 
    \end{subfigure}
    \begin{subfigure}[t]{0.45\textwidth}
    \centering
        \includegraphics[width=\linewidth]{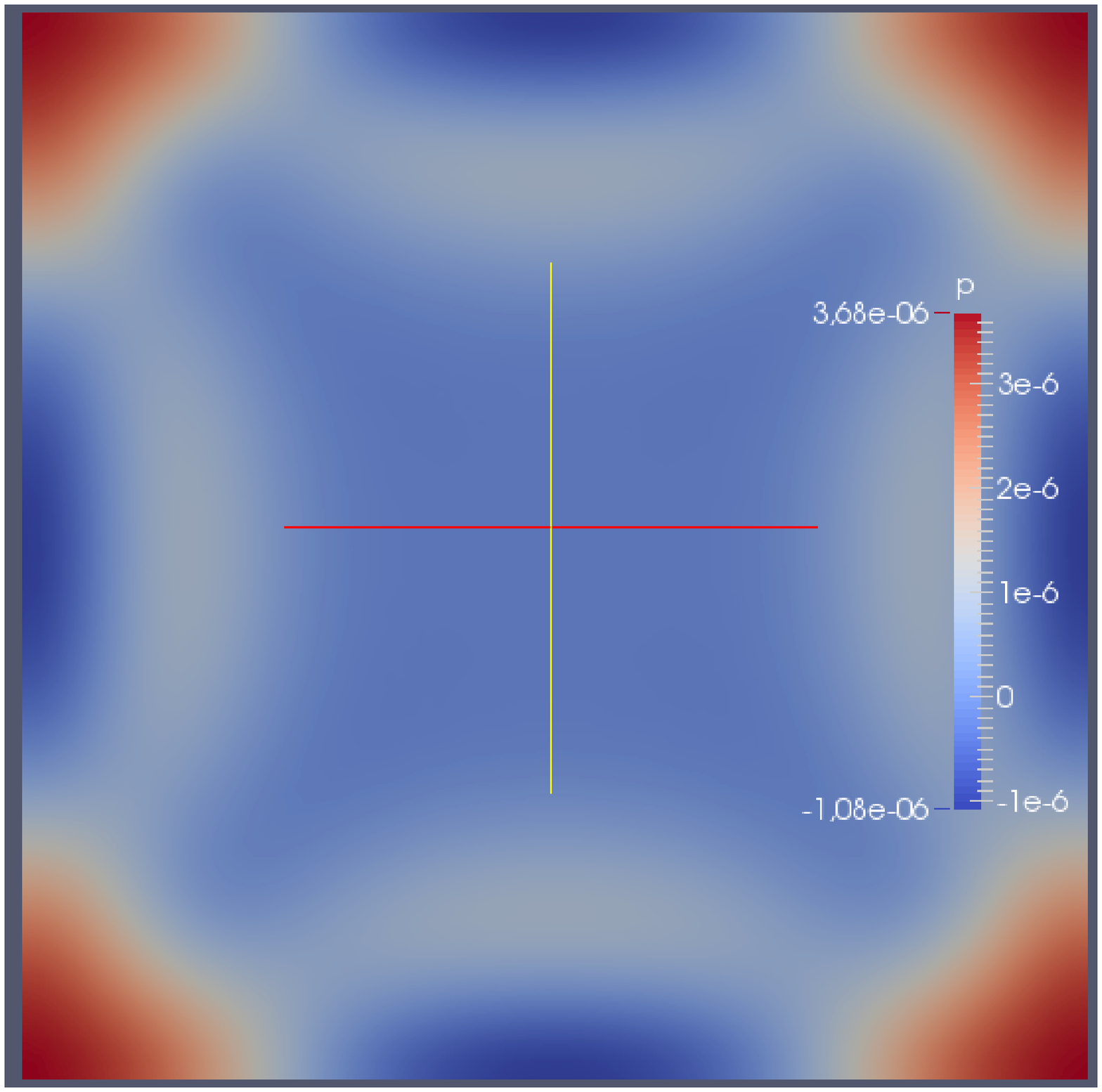}
        \caption{P2 Adomian solution} 
    \end{subfigure}
    \caption{RK-DG Reference and ABS-DG solutions after 6 seconds, for wall
    boundary conditions}
    \label{fig:wallBoundaryAll}
\end{figure}

\begin{figure}[h]
    \centering
    \begin{subfigure}[t]{0.45\textwidth}
        \centering
        \includegraphics[width=\linewidth]{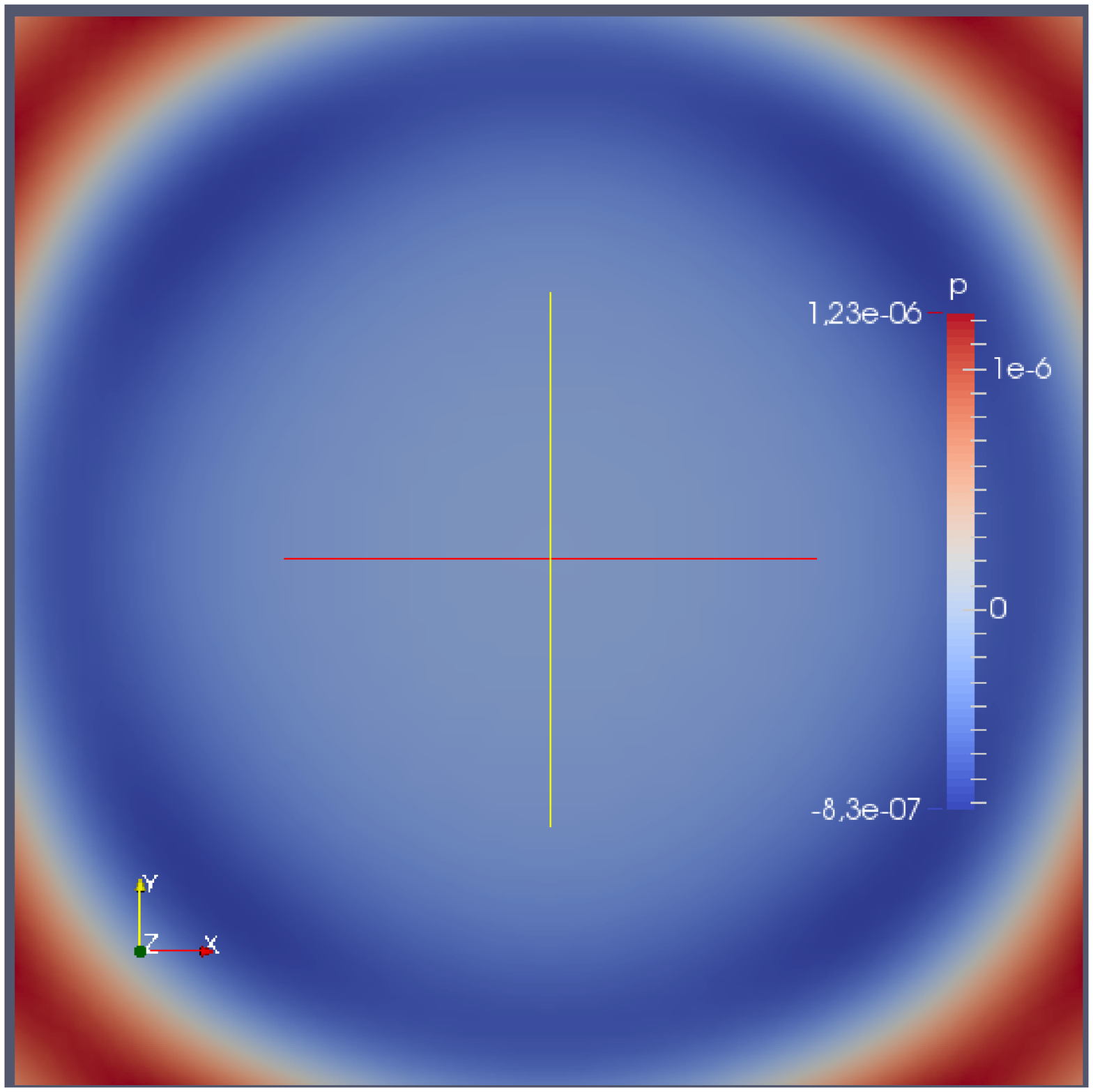}
        \caption{Reference solution}
    \end{subfigure}
    \begin{subfigure}[t]{0.45\textwidth}
        \centering
        \includegraphics[width=\linewidth]{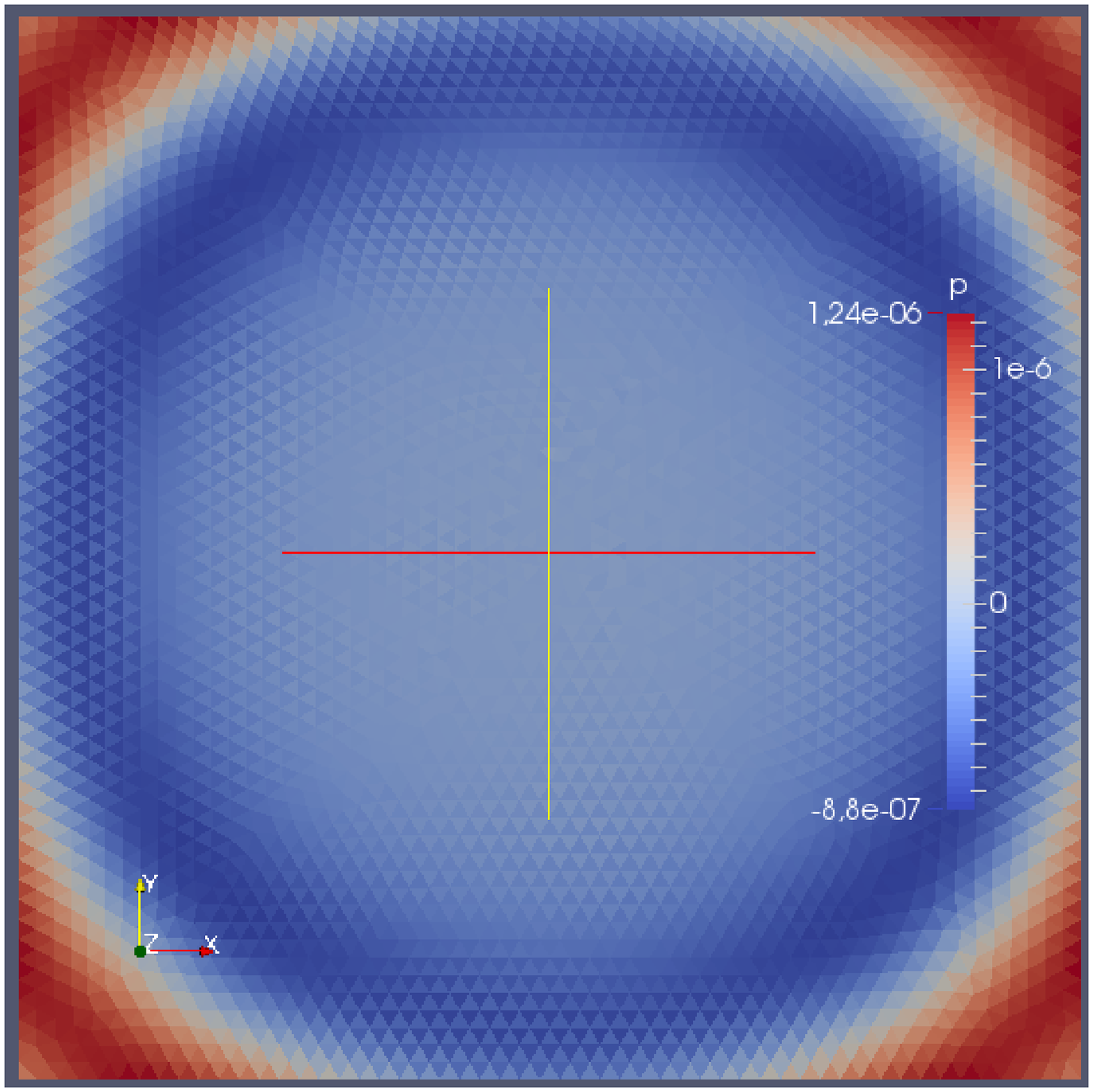}
        \caption{P0 Adomian solution}
    \end{subfigure}
    \begin{subfigure}[t]{0.45\textwidth}
    \centering
        \includegraphics[width=\linewidth]{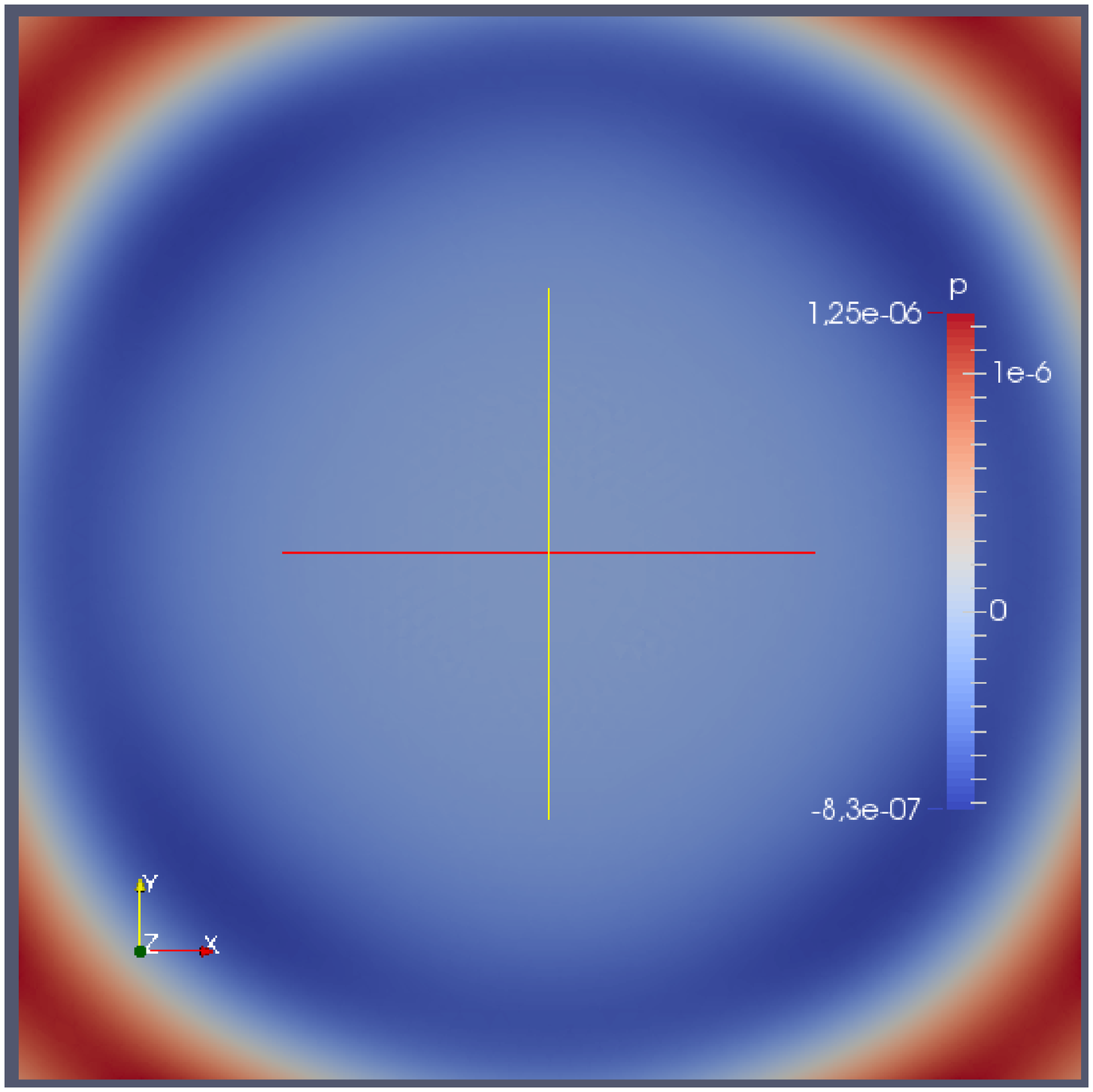}
        \caption{P1 Adomian solution} 
    \end{subfigure}
    \begin{subfigure}[t]{0.45\textwidth}
    \centering
        \includegraphics[width=\linewidth]{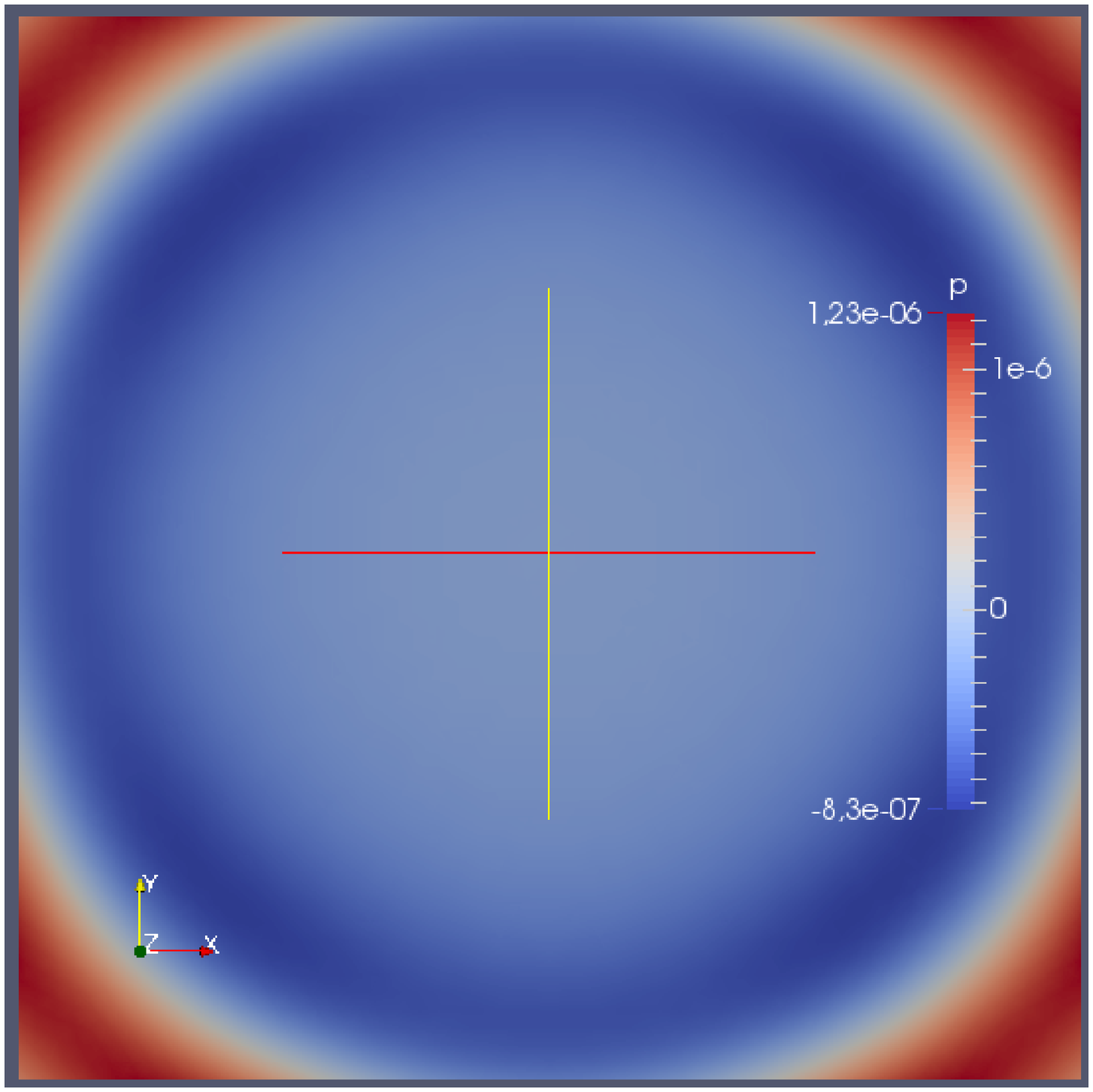}
        \caption{P2 Adomian solution} 
    \end{subfigure}
    \caption{RK-DG Reference and ABS-DG solutions after 6 seconds, for
    non-reflective boundary conditions}
    \label{fig:nonreflectiveBoundaryAll}
\end{figure}

The objective of this test is to show how to appropriately impose boundary
conditions to the ABS-DG method, since the solution is obtained as a series.
Two relevant boundary conditions in aeroacoustic are considered: slip wall and
non-reflective conditions. To make the solution fulfill the imposed boundary
condition, we force each of the ABS-DG series terms to satisfy them. The slip
wall BC are implemented in a weak sense by nullifying the normal to the
boundary component of the momentum flux. Non-reflective BC is achieved
according to the standard characteristic based non-reflecting boundary
conditions. To estimate the accuracy of the ABS-DG and since an exact solution
is not available, a $5^{th}$ order RK-DG simulation is run on a fine grid
(mesh size of $0.1$) and the solution is considered as a reference. The ABS-DG test
is run on a mesh with a size of $0.19$. Simulations are run till $6$ seconds
flow time. Table \ref{BCtest} shows the relative (to reference DG solution) error
for different ABS-DG orders, we can see that we have a very good agreement with
the reference solution and errors are comparable to those obtained in the case
of boundary-free tests. Figures \ref{fig:wallBoundaryAll} and
\ref{fig:nonreflectiveBoundaryAll} show the propagated pulse obtained by the reference
solution and 1st, 2nd and 3rd order (in space) ABS-DG schemes. These tests
demonstrate that imposing the boundary conditions on each term of the series
for the ABS-DG scheme is an appropriate approach.

\begin{table}[ht]
    \caption{ABS-DG results compared to reference DG solution (relative error)} 
    \centering 
    \begin{tabular}{c c c c c} 
        \hline\hline 
        Order & Wall condition & Non-reflective condition \\[0.5ex] 
        \hline 
        1 & 9.39E-002 & 9.97E-002 \\ 
        2 & 3.78E-002 & 2.27E-002 \\
        3 & 3.50E-003 & 1.22E-003 \\[1ex] 
        \hline 
    \end{tabular}
    \label{BCtest} 
\end{table}

%

\newpage
\section{Conclusions}
In this paper a numerical scheme based on the Adomian decomposition method is
proposed (ABS). To assess the method, the space derivative operators are
discretized using the classical discontinuous Galerkin techniques (ABS-DG). The
derivation of the proposed scheme ABS is described in detail, some nice
proprieties are proved making the scheme easy to implement and the integration
in time very accurate. A connection to the Runge-Kutta time discretization
method is established in the linear case, with a clear advantage when using
ABS. Indeed, the proven recursive formula makes ABS (or ABS-DG) easy to
implement as stated above, being the time-order adaptive and dynamic (no need
to set the order in advance) leading to an optimal accuracy with minimum
cost. Finally, the ABS-DG scheme performance is assessed by comparison to the
classical RK-DG results and the exact solution

\section*{Acknowledgments}
This research is supported by the Basque Government through the BERC 2014-2017
program and by the Spanish Ministry of Economy and Competitiveness MINECO: BCAM
Severo Ochoa accreditation SEV-2013-0323.  The authors gratefully acknowledge
the financial support of Diputaci\'{o}n Foral de Bizkaia (DFB) for this
research and the whole BCAM-BALTOGAR project on turbomachinery (grant
BFA/DFB-6/12/TK/2012/00020). Imanol Garcia de Beristain was funded by the
Basque Government Education Department through the Non Doctoral Researcher
Formation Program with reference (PRE\_2013\_1\_1216). Lakhdar Remaki was
partially funded by the Project of the Spanish Ministry of Economy and
Competitiveness with reference MTM2013-40824-P. Alfaisal University grant IRG
with reference IRG16413.

\newpage
\bibliographystyle{asmems4}

\bibliography{lremaki}

\end{document}